\input amstex
\input amsppt.sty
\magnification 1200

\NoBlackBoxes
\pagewidth {13.8cm}
\pageheight {19cm}

\def\sec{\vskip 2mm \noindent}
\def\Ker{\operatorname{Ker}}
\def\Cat{\operatorname{Cat}}
\def\Cub{\operatorname{Cub}}
\def\Tot{\operatorname{Tot}}
\def\Hom{\operatorname{Hom}}
\def\Coker{\operatorname{Coker}}
\def\Grass{\operatorname{\bold {Grass}}}
\def\Sing{\operatorname{Sing}}
\def\Spec{\operatorname{Spec}}
\def\Im{\operatorname{Im}}
\def\Ob{\operatorname{Ob}}
\def\Sk{\operatorname{Sk}}
\def\rk{\operatorname{rk}}

\def\diag{\operatorname{diag}}
\def \ort{\overset \perp \to \oplus}

\def\Tor{\operatorname{Tor}}

\def\exp{\operatorname{exp}}

\def\mod{\operatorname{mod}}

\def\pr{\operatorname{pr}}
\def\ch{\operatorname{ch}}
\def\chw #1{\operatorname{ \widetilde {ch}_{#1}}}
\def\bc{\operatorname{ \widetilde {bc}}}
\def\Tr{\operatorname{Tr}}
\def\tr{\operatorname{tr}}

\def\opartial{{\overline \partial }}
\def\Id {\text{Id}}

\def\Sub {\operatorname {Sub}}
\def\Holim {\operatornamewithlimits {Holim}}
\NoRunningHeads
\baselineskip 15pt
\topmatter
\title
Higher Bott-Chern forms and Beilinson's regulator.  \\
\endtitle
\author    Jose Ignacio Burgos
\footnote{Partially supported by the
DGICYT n$^{o}$. PB93-0790\hfill\break}\\
Steve Wang
\endauthor
\address {J. I. Burgos. Departament d'\` Algebra i Geometria.
Universitat de Bar\-ce\-lo\-na, Gran Via 585. 08007 Barcelona, Spain.
burgos\@cerber.mat.ub.es
}
\endaddress
\address
{S. Wang, 700 E. Henry Clay, Apt. 7 Milwaukee,
Wi 53217, USA.}
\endaddress
\abstract
In this paper, we prove a Gauss-Bonnet theorem for the higher
algebraic $K$-theory of smooth complex algebraic varieties.
To each exact $n$-cube of hermitian vector bundles, we associate a
higher Bott-Chen form, generalizing the Bott-Chern forms associated to
exact sequences. These forms allow us to define characteristic classes
from $K$-theory to absolute Hodge cohomology. Then we prove that these
characteristic classes agree with Beilinson's regulator map.

\endabstract

\endtopmatter
\document

\heading Introduction.
\endheading

The aim of this paper is to generalize the Gauss-Bonnet theorem
to the higher algebraic $K$-theory of smooth complex algebraic varieties.

The Gauss-Bonnet theorem states that, if $X$ is a compact oriented
2-dimensional riemannian manifold, then
$$
\iint_{X} \Cal K\,dA=2\pi \chi (X),
$$
where $\Cal K$ is the gaussian curvature of the riemannian metric and
$\chi (X)$ is the Euler
characteristic of the manifold $X$. Thus, this theorem relates a global
topological invariant, the
Euler characteristic, with a locally defined differential geometrical
object, the gaussian curvature. The Gauss-Bonnet theorem may be restated
saying that the closed form $\Cal K/2\pi $ represents, in de Rham
cohomology, the Euler class of the tangent bundle.

This result was generalized by Chern (\cite {Ch}, see also \cite {Gr},
\cite {M-S}).
If $X$ is an almost complex
manifold and $E$ is a complex vector bundle, by
topological means (for instance obstruction theory), we can define
some chararacteristic classes of the vector bundle, called the Chern
classes of $E$,
$$
c_{j}(E)\in H^{2j}(X,(2\pi i)^{j}\Bbb Z),\qquad j\ge 0.
$$

In this paper we shall restrict the discussion to
the Chern character (see \cite {M-S}), denoted $\ch_{0}(E)$, which
is a certain power series, with rational coefficients, in the Chern
classes. The Chern character
is additive for exact sequences. Nevertheless, since the cohomology
groups we shall consider will be real vector spaces, any result about
the Chern character will imply the analogous result for
any power series in the Chern classes.

Let us provide $E$ with a hermitian metric $h$.
Let $E_{X}$ denote the differential graded commutative algebra of
complex valued
differential forms on $X$, and let $E_{X,\Bbb R}$ denote the subalgebra
of real forms.
Let $D$ be the unique connection of $E$ satisfying
\roster
\item "(1)" $D$ preserves $h$.
\item "(2)" If $U\subset X$ is an open subset and $s$ is a holomorphic
section of $E|_{U}$, then $Ds$ is of pure type $(1,0)$.
\endroster
Let $K=D^{2}$ be the curvature form. Following Chern and Weil, let us
write
$$
\chw{0}(E,h)=\tr \exp (-K)\,\in \bigoplus_{p} (2\pi
i)^{p}E^{p,p}_{X,\Bbb R}. \tag 1
$$
Then the form $\chw{0}(E,h)$ is closed. The Chern theorem states
that the
form $\chw{0}(E,h)$ represents, in de Rham cohomology, the Chern
character class.

Observe that an oriented riemannian 2 dimensional manifold admits a
canonical structure of complex manifold. Since the Euler class is the
top Chern class, the Gauss-Bonnet theorem is a particular case of
Chern's theorem.

The additivity of the Chern character implies that it induces a group
homomorphism
$$
\ch_{0}:K_{0}(X)\longrightarrow \bigoplus_{j}(2\pi i)^{j}H^{2j}(X,\Bbb
Q), $$
where $K_{0}(X)$ is the Grothendiek group of $X$.

Characteristic classes for higher algebraic $K$-theory were
introduced
by Gillet in \cite {Gi}. These classes are defined on any cohomology
theory satisfying certain properties, such as de Rham
cohomology. Nevertheless, in this case, these
higher characteristic classes classes do not give much information. For
instance, for a proper smooth
complex algebraic variety, the only non zero classes are the original
Chern classes from the $K_{0}$ group. In contrast, if the
cohomology theory
is absolute Hodge cohomology, the Chern character map obtained in this
way agrees with the
Beilinson regulator map, which is highly non trivial and is involved in
very deep and far reaching conjectures (\cite {Be}).

Recall that, for $X$ a smooth proper algebraic complex variety, we have
$$
\align
H^{2p}_{\Cal H}(X,\Bbb R(p))&=H^{p,p}(X,\Bbb C)\cap (2\pi
i)^{p}H^{2p}(X,\Bbb R),\\
H^{2p-1}_{\Cal H}(X,\Bbb R(p))&=H^{p-1,p-1}(X,\Bbb C)\cap (2\pi
i)^{p-1}H^{2p-2}(X,\Bbb R).\\
\endalign
$$
Therefore, since, by (1), the Chern character form has the right Hodge
type,
the Chern theorem implies that the Chern character form represents, in
absolute Hodge cohomology, the Chern character class.
The question of generalizing the Chern theorem to the
characteristic classes from higher $K$-theory to absolute Hodge
cohomology arises naturally.

In their paper about arithmetic characteristic classes of hermitian
vector bundles \cite {G-S 1}, Gillet and Soul\'e accomplished the
first step of this program,
extending the
Chern theorem to the case of $K_{1}(X)$. Let us
briefly explain this result.

The elements of $K_{1}(X)$ may be represented by exact sequences of
vector bundles. Thus, the first step is to understand what is the
analogue of the Chern forms of hermitian vector bundles in the case of
exact sequences of hermitian vector bundles.

Let
$$ \overline {\xi }:
0@>>>
(E',h')@>>>
(E,h)@>>>
(E'',h'')@>>>0,
$$
be an exact sequence of hermitian vector bundles.  Then
the Chern character classes satisfy
$$
\ch_{0}(E)=\ch_{0}(E')+\ch_{0}(E'').
$$
Nevertheless, in general, the Chern character form does not behave
additively:
$$
\chw{0}(E,h)\not=\chw{0}(E',h')+\chw{0}(E'',h'').
$$
In the case when $h'$ and $h''$ are the metrics induced by $h$, Bott
and Chern
(\cite {B-C}) have defined a differential form, $\chw{1}(\overline \xi
)$, which will be called the Bott-Chern form of $\overline {\xi }$,
such that
$$
-2\partial \overline \partial \chw{1}(\overline \xi )=
\chw{0}(E',h')+\chw{0}(E'',h'')-\chw{0}(E,h). \tag 2
$$
Note that the normalization factor we use is different from the
normalization factor used in the original paper.
The forms $\chw{1}(\overline \xi )$ are natural and well defined only up
to $\Im \partial +\Im \overline \partial $.

Bismut, Gillet and Soul\'e (\cite {B-G-S}, \cite {G-S 1}) have given a
different construction of Bott-Chern forms that can be applied to the
case when $h'$ and $h''$ are not the induced metrics. These Bott-Chern
forms are also well defined only up to $\Im \partial +\Im \overline
\partial $.

Bott-Chern forms are exactly what we were looking for, and,
when $X$ is proper, Gillet and Soul\'e (\cite
{G-S 1})
have given an explicit description of Beilinson's regulator map for
$K_{1}(X)$ in terms of these Bott-Chern forms, thus extending the Chern
theorem to the group $K_{1}$.

As we have seen,
the Chern
character class on the $K_{0}$ is additive for exact sequences.
Nevertheless one cannot make a consistent choice of
representatives of the Chern character that behave additively for exact
sequences. The Bott-Chern forms measure precisely this lack of
additivity at the level of Chern forms, and they are responsible
for the Chern character for the $K_{1}$-group.
Following Schechtman's ideas (\cite {Sch}) the lack of additivity of the
representatives of the
Chern character for $K_{i}$ is responsible for the Chern character
for $K_{i+1}$.
Thus, the lack of additivity of Bott-Chern forms
should allow us to define second order Bott-Chern forms that give a
description of Beilinson's regulator map for the $K_{2}$. And we can
repeat this process to obtain Beilinson's regulator map for all the $K$
groups.

In this direction, when $X$ is proper, the second author (\cite {Wan})
has defined higher Bott-Chern forms for exact hermitian $n$-cubes. These
forms may be thought of as an iteration of Bott-Chern forms. Moreover,
he has used them to define characteristic classes for higher $K$-theory,
proving that, if one
can naturally extend higher Bott-Chern forms to the non proper case,
then these characteristic classes agree with Beilinson's regulator map.

In this paper we shall give a variant of Wang's original
construction that can be easily extended to
the non-proper case, thus completing the proof of the Gauss-Bonnet
theorem for higher algebraic $K$-theory. An interesting feature of
the construction given here is that we obtain well defined Bott-Chern
forms and not only modulo $\Im \partial +\Im \overline \partial $.

Parallel results in the
framework of multiplicative $K$-theory have been obtained by Karoubi
in \cite {K1} and \cite {K2}.

Gillet and Soul\'e \cite {G-S 1} have used Bott-Chern
forms to define arithmetic $K_{0}$ groups. Soul\'e has suggested
(\cite {So}, see also \cite {De})
that one may define higher arithmetic $K$-groups as the homotopy fibre
of Beilinson's regulator map.
We expect that
higher Bott-Chern forms, as presented in this paper, will be useful in
giving a more concrete definition of higher arithmetic $K$-theory and
studying its properties.

Throughout the paper all vector bundles will be algebraic and we shall
use the equivalent notion of locally free sheaf.

The plan of the paper is as follows. In {\S}1 we recall the definition
of real absolute Hodge cohomology. We shall also show that real absolute
Hodge cohomology can be computed by means of a complex composed by forms
defined on $X\times (\Bbb P^{1})^{n}$, $n\ge 0$. Higher Bott-Chern forms
will live in this complex.

In {\S}2 we introduce and study some properties of smooth at
infinity hermitian metrics. Over a non proper smooth complex variety, to
compute real absolute Hodge cohomology, one needs to impose
logarithmic conditions at infinity to the
differential forms. Thus we cannot
use arbitrary hermitian metrics because they will produce differential
forms with arbitrary singularities at infinity. The use of smooth at
infinity hermitian metrics ensures that Bott-Chern forms
have the right behaviour at infinity.

In {\S}3 we recall the notion of exact metrized $n$-cubes. To each exact
$n$-cube, $E$, we shall attach a vector bundle, $\tr_{n}(E)$ over
$X\times
(\Bbb P^{1})^{n}$ which may be thought of as a homotopy between the
faces
of $E$. When the hermitian metrics of $E$ satisfy certain technical
condition, we shall define a natural metric on $\tr_{n}(E)$. The
Chern character form of the vector bundle $\tr_{n}(E)$ will play the
role of
higher Bott-Chern forms. Note that these forms live in $X\times (\Bbb
P^{1})^{n}$.

In {\S}4 we use higher Bott-Chern forms to define Chern character
classes from higher $K$-theory to real absolute Hodge cohomology.

In {\S}5 we prove that the higher Chern character defined in {\S}4
agrees with Beilinson's regulator map.

In {\S}6 we recall several complexes that compute real absolute Hodge
cohomology and homology. Using them we give, for $X$ proper, two
different versions of higher Bott-Chern forms which are defined on $X$.
The first one, obtained using the Thom-Whitney simple, is
multiplicative.
The
second one agrees with classical Bott-Chern forms and with the original
definition due to Wang.

{\it Acknowledgements.\/} We would like to thank Prof. C. Soul\'e who
suggested
this question to us and helped us with encouragement and numerous hints.
Without him this paper would never have been produced. We would like
to thank
Prof. V. Navarro Aznar for his help and ideas; in particular, the
final definition of Bott-Chern forms is due to a conversation with
him.
Moreover, he proposed some shortcuts in {\S}2.
We would also like to thank
Prof. B. Mazur for his support and guidance and to prof. F. Lecomte
for pointing us some errors in {\S} 5 and the way to fix them.
We acknowledge the
help of many colleagues for useful conversations which have helped us to
understand
a number of aspects of the subject. Our thanks to J.B. Bost, N. Dan, H.
Gillet, D. Grayson, P. Guillen, C. Naranjo, P.Pascual and D. Roessler.

\heading {\S}1 Absolute Hodge cohomology.
\endheading

In this section we shall recall the definition of real absolute Hodge
cohomology \cite {Be} of a smooth complex algebraic variety $X$. By a
smooth
complex variety we shall mean a smooth separated scheme of finite type
over $\Bbb C$. We shall also construct a complex, composed by forms on
$X\times (\Bbb P^{1})^{n}$, $n\ge 0$, whose cohomology is the real
absolute Hodge cohomology of $X$.

\sec
(1.1) Let $\overline X$ be a smooth proper complex variety.
Let $Y\subset \overline X$ be a normal crossing divisor and let us
write $X=\overline X-Y$. Let $E^{*}_{\overline X}$ be the differential
graded algebra of differential forms on $\overline X$, and let
$E_{\overline X}^{*}(\log Y)$ be the
differential graded algebra of C$^{\infty }$ complex differential forms
on $\overline X$
with logarithmic singularities along $Y$ (see \cite {Bu 1}). The algebra
$E_{\overline X}^{*}(\log Y)$ has a real structure,
$E_{\overline X}^{*}(\log Y)_{\Bbb R}$, a weight filtration $W$ defined
over $\Bbb R$
and a Hodge filtration $F$. Moreover the cohomology of this algebra
gives us the cohomology of $X$ with its real mixed Hodge structure.

Let us denote by $\widehat W$ the d\'ecal\'ee
filtration of $W$. That is
$$
\widehat W_{r}E^{n}_{\overline X}(\log Y)=
\{x\in W_{r-n}E^{n}_{\overline X}(\log Y)\mid
dx\in W_{r-n-1}E^{n+1}_{\overline X}(\log Y)\}.
$$

We write
$$
E_{\log}^{*}(X)=\lim\Sb \longrightarrow \\ (\widetilde X_{\alpha },
Y_{\alpha })\endSb E^{*}_{\widetilde X_{\alpha }}(\log Y_{\alpha }),
$$
where the limit is taken along all the smooth compactifications
$\widetilde X_{\alpha }$ of $X$ with $Y_{\alpha }=\widetilde X_{\alpha
}-X$ a normal crossing divisor. Then $E_{\log}^{*}(X)$ is a differential
graded algebra and it has an induced real structure, a weight filtration
and a Hodge filtration. Moreover the map
$$
(E^{*}_{\overline X}(\log Y)_{\Bbb R},\widehat W)
\longrightarrow
(E^{*}_{\log}(X)_{\Bbb R},\widehat W)
$$
is a filtered quasi-isomorphism and the map
$$
(E^{*}_{\overline X}(\log Y),\widehat W,F)
\longrightarrow
(E^{*}_{\log}(X),\widehat W, F)
$$
is a bifiltered quasi-isomorphism.

\sec
(1.2)
Let us write
$$
\frak H^{*}(X,p)=s((2\pi i)^{p}\widehat W_{2p}E^{*}_{\log}(X)_{\Bbb R}
\oplus \widehat W_{2p}\cap F^{p}E^{*}_{\log}(X)@>u>>
\widehat W_{2p}E^{*}_{\log}(X)
),
$$
where $u(r,f)=f-r$ and $s$ denotes the simple of a morphism of
complexes, i.e. the c\^one shifted by one. The differential of this
complex will be denoted by $d_{\frak H}$.

The real absolute Hodge cohomology of $X$ (\cite {Be}) is
$$
H^{n}_{\Cal H}(X,\Bbb R(p))=H^{n}(\frak H(X,p)).
$$

\sec
(1.3) A cubical or cocubical object (see \cite {G-N-P-P}) is an object
modeled on the cube
in the same way as a simplicial or cosimplicial object is modeled on the
simplex. Let $(\Bbb P^{1}_{\Bbb C})^{\cdot}$ be the cocubical scheme
which in degree $n$ is
$(\Bbb P^{1}_{\Bbb C})^{n}$, the $n$-fold product of the complex
projective line. The faces and degeneracies
$$
\alignat 2
d^{i}_{j}:
(\Bbb P^{1}_{\Bbb C})^{n}&\longrightarrow (\Bbb P^{1}_{\Bbb C})^{n+1},\
& i&=1,\dots ,n+1,\ j=0,1\\
s^{i}:
(\Bbb P^{1}_{\Bbb C})^{n}&\longrightarrow (\Bbb P^{1}_{\Bbb C})^{n-1},\
& i&=1,\dots ,n,\\ \endalignat
$$
are given by
$$
\align
d^{i}_{0}(x_{1},\dots ,x_{n})&=
(x_{1},\dots ,x_{i-1},(0:1),x_{i},\dots ,x_{n})\\
d^{i}_{1}(x_{1},\dots ,x_{n})&=
(x_{1},\dots ,x_{i-1},(1:0),x_{i},\dots ,x_{n})\\
s^{i}(x_{1},\dots ,x_{n})&=
(x_{1},\dots ,x_{i-1},x_{i+1},\dots ,x_{n}).\\
\endalign
$$

\sec
(1.4)
The complexes $\frak H^{*}(X\times (\Bbb P^{1})^{\cdot},p)$ form a
cubical complex. We shall write
$$
\align
s_{i}&=(\Id \times s^{i})^{*},\\
d_{i}^{j}&=(\Id \times d^{i}_{j})^{*}.
\endalign
$$
Let us denote by
$
\frak H_{\Bbb P}^{*,*}(X,p)
$
the associated double complex. That is
$$
\frak H_{\Bbb P}^{r,n}(X,p)=
\frak H^{r}(X\times (\Bbb P^{1})^{-n},p),
$$
with differentials $$
\align
d'&=d_{\frak H},\\
d''&=\sum (-1)^{i+j} d^{j}_{i}.
\endalign
$$

\sec
(1.5)
We want to obtain from $\frak H_{\Bbb P}^{*,*}(X,p)$, a complex which
computes the absolute Hodge cohomology of $X$. On the one hand, since we
are
using a cubical theory we need to factor out by the degenerate elements
(see \cite {Mas}). On the other hand, we need to kill all cohomology
classes coming from the projective spaces.

Let us denote by
$p_{0}:X\times (\Bbb P^{1})^{n}\longrightarrow X$
the projection over the first factor and by
$p_{i}:X\times (\Bbb P^{1})^{n}\longrightarrow \Bbb
P^{1}$, $i=1,\dots,n$, the projection over the $i$-th projective
line.

Let $(x:y)$ be homogeneous coordinates of $\Bbb P^{1}$. Let us write
$$
g=\log \frac{x\overline x+y \overline y}{x\overline x}.
$$
Let $\omega =\partial \overline \partial g\in (2\pi i)E^{2}_{\Bbb
P^{1},\Bbb R}$ be a K\"ahler form over $\Bbb
P^{1}$. Let $\omega _{i}=p_{i}^{*}\omega \in E^{*}_{\log}(X\times (\Bbb
P^{1})^{n})$. For an element
$$
x=(r,f,\eta )\in \frak H^{r}(X\times (\Bbb P^{1})^{n},p),
$$
we shall write
$$
\align
\omega _{i}\land x&=(\omega _{i}\land r,\omega _{i}\land f,\omega
_{i}\land \eta)\\
&\in\frak H^{r+2}(X\times (\Bbb P^{1})^{n},p+1).
\endalign
$$

\sec
{\bf Definition 1.1.} We shall denote by
$\widetilde {\frak H}^{*,*}(X,p)$ the double complex given by
$$
\widetilde {\frak H}^{r,n}(X,p)=
\frak H^{r,n}_{\Bbb P}(X,p)\left /
\sum_{i=1}^{-n}s_{i}\left(\frak H^{r,n+1}_{\Bbb P}(X,p)\right)\,\oplus\,
\omega _{i}\land s_{i}\left(\frak H^{r-2,n+1}_{\Bbb P}(X,p-1)\right).
\right.
$$
We shall denote by
$\widetilde {\frak H}^{*}(X,p)$ the associated simple complex. The
differential of this complex will be denoted by $d$.

\sec
In the definition of $\widetilde {\frak H}^{r,n}(X,p)$, the first
summand of the quotient is meant to kill the degenerate
classes, whereas the second summand kill the classes coming
from the projective spaces. The next result shows that we have reached
our objective.

\proclaim{Proposition 1.2}
The natural morphism of complexes
$$
\iota :\frak H^{*}(X,p)=\widetilde {\frak H}^{*,0}(X,p)\longrightarrow
\widetilde {\frak H}^{*}(X,p)
$$
is a quasi-isomorphism.
\endproclaim
\demo{Proof} Since $\widetilde {\frak H}^{*}(X,p)$ is a simple complex
associated to a double complex, there is a second quadrant spectral
sequence with $E_{1}$ term
$$
E_{1}^{n,r}=H^{r}(\widetilde {\frak H}^{*,n}(X,p)).
$$
When this spectral sequence converges, the limit is
$H^{*}(\widetilde {\frak H}^{*}(X,p))$.
The following lemma shows that this spectral sequence converges
and implies that $\iota $ is a quasi-isomorphism.

\proclaim{Lemma 1.3} For $n<0$ the cohomology of the complex
$\widetilde {\frak H}^{*,n}(X,p)$ is zero.
\endproclaim
\demo{Proof} For each $j$ let us write
$$
\widetilde {\frak H}_{j}^{r,n}(X,p)=
\frak H^{r,n}_{\Bbb P}(X,p)\left /
\sum_{i=1}^{j}s_{i}\left(\frak H^{r,n+1}_{\Bbb P}(X,p)\right)\,\oplus\,
\omega _{i}\land s_{i}\left(\frak H^{r-2,n+1}_{\Bbb P}(X,p-1)\right).
\right.
$$
Let us prove, by induction over $j$, that for $j\ge 1$
$$
H^{*}(\widetilde {\frak H}_{j}^{*,n}(X,p))=0.
$$
For $j=1$, $n\le -1$, the complex $\widetilde {\frak
H}_{1}^{r,n}(X,p)$ is the cokernel of the monomorphism
$$
\matrix
\frak H^{*}(X\times (\Bbb P^{1})^{-n-1},p)&
\oplus
&
\frak H^{*}(X\times (\Bbb P^{1})^{-n-1},p-1)[-2]
&\longrightarrow
&\frak H^{*}(X\times (\Bbb P^{1})^{-n},p)\\
\alpha & \oplus &\beta&\longmapsto& s_{1}(\alpha )+\omega _{1}\land
s_{1}(\beta )
\endmatrix
$$
But by the Dold-Thom isomorphism for absolute Hodge cohomology, the
above morphism is a quasi-isomorphism. For $j>1$,
$n<-1$, $\widetilde {\frak H}_{j}^{*,n}(X,p)$ is the cokernel of the
monomorphism
$$
\matrix
\widetilde {\frak H}_{j-1}^{*,n-1}(X,p)
&\oplus &
\widetilde {\frak H}_{j-1}^{*,n-1}(X,p-1)[-2]
&\longrightarrow&
\widetilde {\frak H}_{j-1}^{*,n}(X,p)\\
\alpha &\oplus &\beta
&\longmapsto&
s_{j}(  \alpha )\omega _{j}\land s_{j}(\beta ).
\endmatrix
$$
By induction hypothesis, the source and the target of this morphism have
zero cohomology. Therefore the cokernel also has zero cohomology.

\enddemo

\enddemo

\heading {\S}2 Smooth at infinity hermitian metrics.
\endheading

In this section we introduce smooth at infinity hermitian metrics. For a
smooth complex variety $X$ and a locally free sheaf $\Cal F$, a smooth
at infinity hermitian metric is a metric that can be extended to a
smooth metric over some compactification of $\Cal F$. The interest of
smooth at infinity hermitian metrics is that they provide
representatives of Chern classes in absolute Hodge cohomology.

\sec (2.1)
Before defining smooth at infinity hermitian metrics, we shall study
classes of compactifications of locally free sheaves.

\sec
{\bf Definition 2.1.} Let $X$ be a smooth complex variety
and let $\Cal F$ be a locally free sheaf over $X$. A
compactification of $\Cal F$ is a smooth
compactification of $X$, $i:X\longrightarrow \widetilde X$, a locally
free sheaf
$\widetilde {\Cal F}$ over $\widetilde X$ and an isomorphism $\varphi
:\Cal F \longrightarrow i^{*}\widetilde {\Cal F}$.

A compactification of $\Cal F$ will be denoted by $(i,\widetilde
X,\widetilde {\Cal F}, \varphi )$. Usually, we shall identify $X$ with
$i(X)$ and $\Cal F$ with $\widetilde {\Cal F}|_{X}$, and denote a
compactification by $(\widetilde {\Cal F},\widetilde X)$.

\proclaim{Proposition 2.2} Let $X$ be a smooth complex variety
and let $\Cal F$ be a locally free sheaf over $X$. Then there
exists a compactification of $\Cal F$. \endproclaim
\demo{Proof} Let $X\longrightarrow \widetilde X_{1}$ be any
compactification of $X$. Then there is a coherent sheaf $\widetilde
{\Cal F}_{1}$ on $\widetilde X_{1}$ such that $\widetilde {\Cal
F}_{1}|_{X}=\Cal F$. By \cite {Ro} (see also \cite {Ri} and \cite
{N 1}) there is a proper modification $\psi :\widetilde X\longrightarrow
\widetilde X_{1} $, which induces an isomorphism $\psi
^{-1}(X)\longrightarrow X$, and such that $\widetilde {\Cal F}=\psi
^{*}(\widetilde
{\Cal F}_{1})/\Tor (\psi ^{*}(\widetilde {\Cal F}_{1}))$ is a locally
free sheaf. Moreover $\widetilde {\Cal F}|_{\psi ^{-1}(X)}$ is
isomorphic to $\widetilde {\Cal F}_{1}|_{X}$. Thus the induced map
$i:X\longrightarrow
\widetilde X$ is a compactification of $X$, and $\widetilde {\Cal F}$ is a
compactification of $\Cal F$.
\enddemo

\sec
{\bf Definition 2.3.} Let $X$ be a smooth complex variety
and let $\Cal F$ be a locally free sheaf over $X$. Let
$(i_{1},\widetilde X_{1},\widetilde {\Cal F}_{1},\varphi _{1})$ and
$(i_{2},\widetilde X_{2},\widetilde {\Cal F}_{2},\varphi _{2})$ be two
compactifications of $\Cal F$. We say that $\widetilde {\Cal F}_{1}$ and
$\widetilde {\Cal F}_{2}$ are equivalent if there exists a third
compactification
$(i_{3},\widetilde X_{3},\widetilde {\Cal F}_{3},\varphi _{3})$ and
morphisms $\psi _{1}:\widetilde X_{3}\longrightarrow \widetilde X_{1}$
and $\psi _{2}:\widetilde X_{3}\longrightarrow \widetilde X_{2}$ such
that
\item {1)} $\psi_{1} \circ i_{3}=i_{1}$ and  $\psi_{2} \circ
i_{3}=i_{2}$.
\item {2)} There are isomorphisms $\alpha _{1}:\widetilde {\Cal F}_{3}
\longrightarrow \psi _{1}^{*}\widetilde {\Cal F}_{1}$ and
$\alpha _{2}:\widetilde {\Cal F}_{3}
\longrightarrow \psi _{2}^{*}\widetilde {\Cal F}_{2}$ such that
$i_{3}^{*}\alpha _{1} \circ \varphi _{3}=\varphi _{1}$ and
$i_{3}^{*}\alpha _{2} \circ \varphi _{3}=\varphi _{2}$.

In order to simplify the notation, a class of equivalent
compactifications of $\Cal F$ will be denoted by a single symbol, for
instance $\widetilde {\Cal F}$.
Moreover, if there is no danger of confusion, we shall denote by the
same
symbol the locally free sheaf which appears in any representative of
this class.

\sec
(2.2) Let us see that a compactification class induces uniquely
determined compactification classes in quotients and subsheaves.

\proclaim{Theorem 2.4} Let $X$ be a smooth complex variety
and let
$$
\xi :\  0@>>> \Cal F @>f>>\Cal G @>g>> \Cal H@>>>0
$$
be an exact sequence of locally free sheaves over $X$. Then, for any
compactification class $\widetilde {\Cal G}$ of $\Cal G$, there are
uniquely determined compactification classes $\widetilde {\Cal F}$
and $\widetilde {\Cal H}$ of $\Cal F$ and $\Cal H$ respectively, such
that $\xi $ extends to an exact sequence
$$
\widetilde {\xi} :\  0@>>> \widetilde{\Cal F}
@>\widetilde f>>\widetilde {\Cal G}
@>\widetilde g>> \widetilde {\Cal H}@>>>0,
$$
over a compactification $\widetilde X$ of $X$.
\endproclaim
\demo{Proof}
Let $\widetilde X_{1}$ be a compactification of $X$ where $\widetilde
{\Cal G}$ is defined. Let $r=\rk \Cal H$. Let
$\Grass_{\widetilde X_{1}}^{r}(\widetilde {\Cal G})$ be the Grassmanian
of rank
$r$ quotients of $\widetilde {\Cal G}$ (\cite {G-D}). Let us denote by
$\Cal U$ the universal bundle on
$\Grass_{\widetilde X_{1}}^{r}(\widetilde {\Cal G})$. The
exact sequence
$$
\xi :\  0@>>> \Cal F @>f>>\Cal G @>g>> \Cal H@>>>0
$$
induces a morphism
$$
\varphi :X\longrightarrow
\Grass_{\widetilde X_{1}}^{r}(\widetilde {\Cal G}).
$$
By resolution of singularities,
there is a proper modification $\widetilde X$ of $\widetilde X_{1}$,
which is a compactification of $X$ and such that $\varphi $ extends to a
morphism
$$
\widetilde \varphi :\widetilde X\longrightarrow
\Grass_{\widetilde X_{1}}^{r}(\widetilde {\Cal G}).
$$
Then $\widetilde {\Cal H}=\widetilde \varphi ^{*}(\Cal U)$ is a
compactification of $\Cal H$, $\widetilde {\Cal F}=\Ker(\widetilde
{\Cal G}\longrightarrow \widetilde {\Cal H})$ is a compactification of
$\Cal F$ and $\xi $ extends to an exact sequence
$$
\widetilde {\xi} :\  0@>>> \widetilde{\Cal F}
@>\widetilde f>>\widetilde {\Cal G}
@>\widetilde g>> \widetilde {\Cal H}@>>>0.
$$

The unicity follows from the fact that, since $X$ is dense in
$\widetilde X$, the morphism $\widetilde \varphi $ is unique.

\enddemo

\sec
{\bf Definition 2.5.} Let $X$ be a smooth complex variety
and let
$$
\xi :\  0@>>> \Cal F @>f>>\Cal G @>g>> \Cal H@>>>0
$$
be an exact sequence of locally free sheaves over $X$. Let $\widetilde
{\Cal G}$ be a class of compactifications of $\Cal G$. Then the
classes of compactifications $\widetilde {\Cal F}$
and $\widetilde {\Cal H}$, of $\Cal F$ and $\Cal H$ respectively,
obtained in theorem 2.4 are called the induced compactifications.

\sec
(2.3) Let us introduce smooth at infinity hermitian metrics.
\sec
{\bf Definition 2.6.} Let $X$ be a smooth complex variety,
let $\Cal F$ be a locally free sheaf over $X$ and let $h$ be
an hermitian metric on $\Cal F$. We say that
$h$ is smooth at infinity if there exist a compactification $\widetilde
{\Cal F}$ of $\Cal F$, and a smooth metric $\widetilde h$ on $\widetilde
{\Cal F}$ such that $\widetilde h|_{X}=h$.

A smooth at infinity hermitian metric determines univocally a
compactification class.

\proclaim{Proposition 2.7}
Let $X$ be a smooth complex variety
and let $\Cal F$ be a locally free sheaf on $X$.
Let $\widetilde {\Cal F}$ and $\widetilde {\Cal F}'$ be two
compactifications of $\Cal F$ and let $\widetilde h$ and $\widetilde
h'$ be smooth metrics on $\widetilde {\Cal F}$ and $\widetilde {\Cal
F}'$.  If
$\widetilde h|_{X}=\widetilde h'|_{X}$,
then
$\widetilde {\Cal F}$ and $\widetilde {\Cal F}'$ are equivalent
compactifications.
\endproclaim
\demo{Proof} We can assume that both compactifications are defined over
the same variety $\widetilde X$. Let $\Cal K_{\widetilde X}$  be the
sheaf of rational functions over $\widetilde X$.

The identity on $\Cal F$ induces morphisms
$$
\align
f:
\widetilde {\Cal F}\otimes \Cal K_{\widetilde X} &\longrightarrow
\widetilde {\Cal F}'\otimes \Cal K_{\widetilde X},\\
f':
\widetilde {\Cal F}'\otimes \Cal K_{\widetilde X} &\longrightarrow
\widetilde {\Cal F}\otimes \Cal K_{\widetilde X},
\endalign
$$
which are inverses of each other. By symmetry it is enough to show
that $f(\widetilde {\Cal F})\subset \widetilde {\Cal F}'$.

Let $U$ be a
Zariski open subset of $\widetilde X$. A section $s\in \Gamma (
U,\widetilde {\Cal F}'\otimes \Cal K_{\widetilde X})$ belongs to
$\Gamma (U,\widetilde {\Cal F}')$ if and only if $\widetilde
h'(s(x))<\infty $ for all $x\in U$. But if $s\in \Gamma (U,\widetilde
{\Cal
F})$ then $\widetilde h'(f(s))|_{X\cap U}=\widetilde h(s)|_{X\cap U}$.
Since $U\cap X$ is dense in $U$ we have
$\widetilde h'(f(s(x)))=\widetilde h(s(x))<\infty $ for all $x\in U$.
\enddemo

\proclaim{Proposition 2.8} Let
$$
\xi : 0@>>>\Cal F'@>>> \Cal F @>>> \Cal F'' @>>> 0
$$
be an exact sequence of locally free sheaves on $X$ and let $h$ be a
smooth at infinity metric on $\Cal F$. Then the metrics $h'$ and
$h''$ induced by $h$ in $\Cal F'$ and $\Cal F''$ are smooth at infinity.
\endproclaim
\demo{Proof} Let $\widetilde {\Cal F}$ be a compactification of $\Cal F$
provided with a metric $\widetilde h$, such that $\widetilde h|_{X}=h$.
By theorem 2.4. there are compactifications $\widetilde {\Cal F}'$ and
$\widetilde {\Cal F}''$ such that $\xi $ can be extended to an exact
sequence
$$
\widetilde \xi :
0@>>>\widetilde {\Cal F}'@>>> \widetilde {\Cal F}
@>>> \widetilde {\Cal F}'' @>>> 0.
$$
Then the metric $\widetilde h$ induces smooth metrics $\widetilde h'$
and $\widetilde h''$ on $\widetilde {\Cal F}'$ and $\widetilde {\Cal
F}''$.
But the restrictions of $\widetilde h'$ and $\widetilde h''$ to $X$ are
$h'$ and $h''$. Therefore these metrics are smooth at infinity.
\enddemo

\proclaim{Proposition 2.9} Let $f:X\longrightarrow Y$ be a morphism
between smooth complex varieties. Let $(\Cal F,h)$ be a locally
free sheaf over $Y$ with $h$ a smooth at infinity metric. Then
$(f^{*}h)$ is a smooth at infinity metric on the locally free
sheaf $f^{*}\Cal F$. \endproclaim
\demo{Proof} Let $(\widetilde Y,\widetilde {\Cal F})$ be a
compactification  of $(Y,\Cal F)$, such that there exists a hermitian
metric $\widetilde h$ with $\widetilde h|_{X}=h$. Let $\widetilde
X_{1}$ be any compactification of $X$. We shall denote by $\Gamma $
the graph of $f$, and by
$\overline \Gamma $ the adherence of $\Gamma $ in $\widetilde
X_{1}\times \widetilde Y$. Let $\widetilde X$ be a resolution of
singularities of $\overline \Gamma $ and let $\widetilde f:\widetilde
X\longrightarrow \widetilde Y$ be the induced morphism. Then
$(\widetilde X,\widetilde f^{*}\widetilde {\Cal F})$ is a
compactification of $(X,f^{*}\Cal F)$ and $\widetilde f^{*}\widetilde h$
is a smooth metric such that $\widetilde f^{*}\widetilde h|_{X}=f^{*}h$.
Therefore $f^{*}h$ is smooth at infinity.

\enddemo

\sec
(2.4) Let us see that smooth at infinity hermitian metrics provide
representatives of the Chern character classes in absolute Hodge
cohomology. Let $X$ be a smooth complex variety, $\Cal F$ a locally free
sheaf and $h$ a smooth at infinity hermitian metric. Let $\widetilde
{\Cal F}$ be the compactification class of $\Cal F$ determined by $h$,
$\widetilde X$ a compactification of $X$ where $\widetilde {\Cal F}$ is
defined, and $\widetilde h$ a smooth metric on $\widetilde {\Cal F}$
extending $h$. Let $K$ (resp. $\widetilde K$) be the curvature form of
$(\Cal F,h)$ (resp. $(\widetilde {\Cal F},\widetilde h)$). Let us write
$$
\align
\chw{0}(\Cal F,h)&=\Tr \exp (-K),\\
\chw{0}(\widetilde {\Cal F},\widetilde h)&=\Tr \exp (-\widetilde K).
\endalign
$$
These forms are closed. Moreover,
$$
\chw{0}(\widetilde {\Cal F},\widetilde h)\in
\bigoplus (2\pi i)^{p} E_{\widetilde X,\Bbb R}^{p,p}.
$$
Since $\chw{0}(\widetilde {\Cal F},\widetilde h)|_{X}=
\chw{0}(\Cal F,h)$,
$$
\chw{0}(\Cal F,h)
\in \bigoplus_{p\ge 0} \left(
W_{0} E^{2p}_{\log}(X)_{\Bbb R}\cap
W_{0}\cap F^{p} E^{2p}_{\log}(X)
\right).
$$
Since this form is closed,
$$
\chw{0}(\Cal F,h)
\in \bigoplus_{p\ge 0} \left(
\widehat W_{2p} E^{2p}_{\log}(X)_{\Bbb R}\cap
\widehat W_{2p}\cap F^{p} E^{2p}_{\log}(X)
\right).
$$
Thus the triple
$$
\chw{0}(\Cal F,h)_{\Cal H}=
(\chw{0}(\Cal F,h),
\chw{0}(\Cal F,h),0)
$$
is a cycle of $ \bigoplus_{p\ge 0} \frak H^{2p}(X,p) $.

\proclaim{Proposition 2.10} The cycle $\chw{0}(\Cal F,h)_{\Cal H}$
represents the Chern character of $\Cal F$ in absolute Hodge cohomology.
\endproclaim
\demo{Proof} If $X$ is proper we have
$$
H^{2p}_{\Cal H}(X,\Bbb R(p)) =H^{p,p}(X,(2\pi i)^{p}\Bbb R).
$$
Therefore the result follows from the classical description of the Chern
character in terms of curvature forms. In the non proper case it follows
from the functoriality of the Chern character.
\enddemo

\heading {\S}3 Exact $n$-cubes of locally free sheaves.
\endheading
In this section we shall recall the notion of exact $n$-cube (see \cite
{Lo 2}, \cite {Wan}). To each metrized exact $n$-cube, $\overline {\Cal
F }$, which satisfies certain conditions, we shall
associate a metrized locally free sheaf on $X\times (\Bbb P^{1})^{n}$,
called the $n$-th transgression of $\overline {\Cal F}$. This
transgression can be viewed as a homotopy between its vertexes. The
Chern character form of the transgression will play the role of higher
Bott-Chern forms.

\sec
(3.1)
First some notations. Let $\left< -1,0,1 \right>$ be the category
associated to the ordered set $\{-1,0,1\}$. Let
$\left< -1,0,1 \right>^{n}$ be its $n$-th cartesian power. By
convention, the
category $\left< -1,0,1 \right>^{0}$ has one element and one morphism.

Let $\frak E$ be an exact category.
\sec
{\bf Definition 3.1.} A $n$-cube of $\frak E$, $\Cal F$,
is a functor from
$\left< -1,0,1 \right>^{n}$ to $\frak E$.

\sec
{\bf Definition 3.2.} Given a $n$-cube $\Cal F$, and numbers $i\in
\{1,\dots ,n\}$, $j\in \{-1,0,1\}$, then the $n-1$-cube,
$\partial ^{j}_{i}\Cal F$ defined by
$$
(\partial ^{j}_{i}\Cal F)_{\alpha _{1},\dots ,\alpha _{n-1}}=
\Cal F_{\alpha _{1},\dots ,\alpha _{i-1},j,\alpha _{i},\dots ,\alpha
_{n-1}} $$
is called a face of $\Cal F$. Given a number $i\in \{1,\dots ,n\}$ and a
$n-1$-tuple
$\alpha =(\alpha _{1},\dots ,\alpha _{n-1})\in \{ -1,0,1 \}^{n-1}$,
the sequence
$$
\partial ^{\alpha }_{i^{c}}\Cal F=
\partial ^{\alpha _{n-1}}_{n}
\dots
\partial ^{\alpha _{i}}_{i+1}
\partial ^{\alpha _{i-1}}_{i-1}
\dots
\partial ^{\alpha _{1}}_{1} \Cal F
$$
is called an edge of $\Cal F$.

Explicitly, the edge
$\partial ^{\alpha }_{i^{c}}\Cal F$ is
$$
\Cal F_{\alpha _{1},\dots ,\alpha _{i-1},-1,\alpha _{i},\dots ,\alpha
_{n-1}}
@>>>
\Cal F_{\alpha _{1},\dots ,\alpha _{i-1},0,\alpha _{i},\dots ,\alpha
_{n-1}}
@>>>
\Cal F_{\alpha _{1},\dots ,\alpha _{i-1},1,\alpha _{i},\dots ,\alpha
_{n-1}}.
$$

{\bf Definition 3.3.} A $n$-cube is called an exact $n$-cube if all its
edges are short exact sequences.

We shall denote by $\underline
C_{n}\frak E$ the exact category of exact $n$-cubes.
Observe that, for all non negative integers $n$, $m$, there is a natural
isomorphism of categories $\underline
C_{n}\underline C_{m}\frak E\longrightarrow \underline C_{n+m}\frak E$.
In particular, an exact $n$-cube can be viewed as an exact
sequence of exact $n-1$-cubes or as an exact $n-1$-cube of exact
sequences.

The maps
$$
\partial ^{j}_{i}
:\Ob \underline C_{n}\frak E\longrightarrow
\Ob\underline C_{n-1}\frak E,
$$
are called face maps.
The maps
$$
s^{j}_{i}
:\Ob \underline C_{n}\frak E\longrightarrow
\Ob\underline C_{n+1}\frak E, \qquad\text{for }i=1,\dots
,n, \text{ and }j=-1,1,
$$
given by
$$
s^{j}_{i}(\Cal F)_{\alpha _{1},\dots ,\alpha _{n+1}}=
\cases
0,&\text{ if } \alpha _{i}=j,\\
\Cal F_{\alpha _{1},\dots ,\alpha _{i-1},\alpha _{i+1},\dots ,\alpha
_{n+1}},&\text{ if }\alpha _{i}\not=j,
\endcases
$$
are called degeneracy maps. An exact $n$-cube $\Cal F\in \Im s
^{j}_{i}$ is called degenerate.

\sec
(3.2) We shall write $C_{n}\frak E=\Ob \underline C_{n}\frak E$ and
$C\frak E=\coprod C_{n}\frak E$.

Assume that the category $\frak E$ is small. To avoid set theoretical
problems, in the sequel we shall always assume tacitly that we replace
any large category by an equivalent small full subcategory.
Observe that the diagram $C\frak E$ behaves like a cubical diagram. We
have replaced the category $\left<0,1\right>$ by the category
$\left<-1,0,1\right>$. This motivates the following construction.

Let $\Bbb ZC_{n}\frak E$ be the
free abelian group generated by $C_{n}\frak E$. And let the differential
$d:\Bbb ZC_{n}\frak E\longrightarrow \Bbb ZC_{n-1}\frak E$ be given by
$$
d=\sum_{i=1}^{n}\sum_{j=-1}^{1} (-1)^{i+j}\partial _{i}^{j}.
$$
Let $D_{n}\subset \Bbb ZC_{n}\frak E$ be the subgroup generated by the
degenerate exact $n$-cubes. Then $dD_{n}\subset D_{n-1}$. Therefore the
following definition makes sense.

\sec
{\bf Definition 3.4.} The homology complex associated to $C\frak E$ is
$$
\widetilde {\Bbb Z}C\frak E=\left. \Bbb ZC\frak E\right/D.
$$

\sec
(3.3)
For the remainder of the section, let us fix a smooth complex variety
$X$. Let $\frak E(X)$
be the exact category of locally free sheaves on $X$ and let $\overline
{\frak E}(X)$ be the exact category of pairs $(\Cal F,h)$, where $\Cal
F\in \Ob \frak E(X)$ and $h$ is a smooth at infinity hermitian metric on
$\Cal F$. The morphisms of this category are
$$
\Hom_{\overline {\frak E}(X)}((\Cal F,h),(\Cal F',h'))=
\Hom_{\frak E(X)}(\Cal F,\Cal F').
$$
Let $F$ be the forgetful functor
$\overline {\frak E}(X)\longrightarrow
\frak
E(X)$.
By choosing metrics we may construct a functor $G:
\frak E(X)\longrightarrow
\overline {\frak
E}(X)$. Then $F\circ G$ is the identity functor on the category $\frak
E(X)$. Moreover, the identity morphism on the vector bundles is a
natural transformation between $G\circ F$ and the identity functor on
$\overline {\frak E}(X)$. Thus $F$ is an equivalence of categories.

For simplicity we shall write $C(X)=C\overline {\frak E}(X)$.
An element
$\overline {\Cal F}\in C_{n}(X)$
is called a metrized exact
$n$-cube of locally free sheaves.

\sec
(3.4)
For technical reasons we need to work with metrized exact $n$-cubes
which have, in all the quotients, the induced metrics.

\sec
{\bf Definition 3.5.} We shall say that a metrized exact $n$-cube,
$\overline {\Cal F}=\{(\Cal F_{\alpha },h_{\alpha })\}$ has induced
quotient metrics (an emi-$n$-cube for short)
if, for each $n$-tuple $\alpha =(\alpha _{1},\dots
,\alpha _{n})$, and each $i$ with $\alpha
_{i}=1$, the metric $h_{\alpha }$ is induced by the metric $h_{(\alpha
_{1},\dots ,\alpha _{i-1},0,\alpha _{i+1},\dots ,\alpha _{n})}$.

Let us see that there are enough emi-$n$-cubes.
Let $\alpha \in \{ -1,0,1 \}^{n}$ be a $n$-tuple. We shall
write $\alpha \le 0$ if $\alpha _{i}\le 0$ for all $i$.

\proclaim{Proposition 3.6}
Let $\Cal F$ be an exact
$n$-cube of locally free sheaves and, for all $\alpha \le 0$, let
$h_{\alpha }$ be a hermitian metric on $\Cal F_{\alpha }$. Then there
is a unique way
to choose metrics $h_{\alpha }$ for all $\alpha \not \le 0$, such that
$\overline
{\Cal F}=\{(\Cal F_{\alpha },h_{\alpha })\}$ is an emi-$n$-cube.
\endproclaim
\demo{Proof} The uniqueness is clear. For the existence, we have to see
that, in each $\Cal F_{\alpha }$, with $\alpha \not \le 0$, all the
possible induced metrics agree. This is guaranteed by the following
result.

\proclaim{Lemma 3.7} Let $\{E_{i,j}\}_{i,j=-1,0,1}$ be an exact
$2$-cube of complex vector spaces. Let $h$ be a hermitian metric on
$E_{0,0}$ and let $h_{1,0}$ and $h_{0,1}$ be the hermitian
metrics in $E_{1,0}$ and $E_{0,1}$ induced by $h$. Then the metrics
induced by $h_{1,0}$ and $h_{0,1}$ in $E_{1,1}$ agree.
\endproclaim
\demo{Proof} Let us identify $E_{-1,0}$ and $E_{0,-1}$ with their images
in $E_{0,0}$. Then the metric $h_{1,0}$ in $E_{1,0}$ is
induced by the isomorphism $E_{-1,0}^{\perp}\cong E_{1,0}$. Therefore we
can identify $E_{1,0}$ with $E_{-1,0}^{\perp}$ and the morphism
$E_{0,0}\longrightarrow E_{1,0}$ with the orthogonal projection. But the
image of $E_{0,-1}$ by this orthogonal projection
is $(E_{-1,0}+E_{0,-1})\cap E_{-1,0}^{\perp}$.
Therefore the metric in $E_{1,1}$ induced by $h_{1,0}$ is
induced by the isomorphism $(E_{-1,0}+E_{0,-1})^{\perp}\cong E_{1,1}$.
By symmetry, the same is true for the metric induced by $h_{0,1}$.
\enddemo
\enddemo

\sec
(3.5)
Let $\Bbb ZC_{emi}(X)$ be the subcomplex of
$\Bbb ZC(X)$ generated by the emi-$n$-cubes, and let
$D_{emi}$ be the subcomplex of $\Bbb ZC_{emi}(X)$ generated by
the degenerate emi-$n$-cubes. We shall write
$$
\widetilde {\Bbb Z}C_{emi}(X)=\left.
\Bbb ZC_{emi}(X)\right/D_{emi}
\subset
\widetilde {\Bbb Z}C(X).
$$

To translate results about emi-$n$-cubes to all exact metrized $n$-cubes
we need to construct a morphism of complexes
$$
\widetilde {\Bbb Z}C(X)\longrightarrow
\widetilde {\Bbb Z}C_{emi}(X).
$$

If $\alpha \in \{-1,0,1\}^{n}$ with $\alpha _{i}>-1$, we shall
write
$\alpha -i=(\alpha _{1},\dots ,\alpha _{i}-1,\dots ,\alpha
_{n} )$.
Let $\overline {\Cal F}=\{(F_{\alpha },h_{\alpha })\}\in C_{n}(X)$.
For $i=1,\dots ,n$ let $\lambda ^{1}_{i}\overline {\Cal F
}$ be defined by
$$
\lambda ^{1}_{i}\overline {\Cal F}_{\alpha }=
\cases
(\Cal F_{\alpha },h_{\alpha }),&\text{ if }\alpha _{i}=-1,0,\\
(\Cal F_{\alpha },h'_{\alpha }),&\text{ if }\alpha _{i}=1,
\endcases
$$
where $h'_{\alpha }$ is the metric induced by $h_{\alpha -i}$. Thus
$\lambda _{i}^{1}\overline {\Cal F}$ has the same locally free sheaves
as $\overline {\Cal F}$, but we have replaced the metrics of the
locally free sheaves of the face $\partial _{i}^{1}\overline {\Cal F}$,
by the metrics induced by $\partial _{i}^{0}\overline {\Cal F}$.

Let $\lambda ^{2}_{i}\overline {\Cal F}$ be the exact
$n$-cube determined by
$$
\align
\partial ^{-1}_{i}\lambda ^{2}_{i}\overline {\Cal F}&=
\partial ^{1}_{i}\overline {\Cal F},\\
\partial ^{0}_{i}\lambda ^{2}_{i}\overline {\Cal F}&=
\partial ^{1}_{i}\lambda ^{1}_{i}\overline {\Cal F},\\
\partial ^{1}_{i}\lambda ^{2}_{i}\overline {\Cal F}&=
0.\\
\endalign
$$
This $n$-cube measures in some sense the difference between $\overline
{\Cal F}$ and $\lambda _{i}^{1}\overline {\Cal F}$.

Let us write $\lambda _{i}\overline {\Cal F}=
\lambda ^{1}_{i}\overline {\Cal F}+\lambda ^{2}_{i}\overline {\Cal F}$,
and let us denote by $\lambda $ the map
$$
\matrix
\lambda :&\Bbb ZC_{n}(X)&\longrightarrow&
\Bbb ZC_{n}(X)\\
&\overline {\Cal F}&\longmapsto&
\cases
\lambda _{n}\dots \lambda _{1}\overline {\Cal F},&\text{ if }n\ge 1,\\
\overline {\Cal F},&\text{ if }n=0.
\endcases
\endmatrix
$$
Then one can check the following properties:
\roster
\item  $\lambda $ is a morphism of complexes.
\item $\Im \lambda \subset \Bbb ZC_{emi}(X)$.
\item $\lambda (D)\subset D_{emi}$.
\endroster
Therefore this map induces a morphism of complexes
$$
\lambda :
\widetilde {\Bbb Z}C(X)
\longrightarrow
\widetilde {\Bbb Z}C_{emi}(X).
$$
In fact $\lambda $ is a homotopy equivalence. The inverse being the
inclusion
$\widetilde {\Bbb Z}C_{emi}(X)
\longrightarrow
\widetilde {\Bbb Z}C(X)$.

\sec
(3.6) Let $\overline {\Cal F}$ be an emi-$n$-cube of locally free
sheaves. We shall associate to it a locally free sheaf
$\tr_{n}(\overline {\Cal F})$
on $X\times (\Bbb P^{1})^{n}$ which, roughly speaking, is a homotopy
between the vertexes of $\overline {\Cal F}$.

Let $((x_{1}:y_{1}),\dots ,(x_{n}:y_{n}))$ be homogeneous coordinates of
$(\Bbb P^{1})^{n}$. Let $\Cal I_{x_{i}}$ (resp. $\Cal I_{y_{i}}$) be the
sheaf of ideals in $X\times (\Bbb P^{1})^{n}$ defined by the subvariety
$x_{i}=0$, (resp. $y_{i}=0$). Let $p_{0}:X\times (\Bbb
P^{1})^{n}\longrightarrow X$ and $p_{i}:X\times (\Bbb
P^{1})^{n}\longrightarrow \Bbb P^{1}$, $i=1,\dots ,n$, be the
projections. Then the maps
$$
\align
\Cal I_{x_{i}}& @>x_{i}^{-1}>> p_{i}^{*}\Cal O_{\Bbb P^{1}}(-1),\\
\Cal I_{y_{i}}&@>y_{i}^{-1}>> p_{i}^{*}\Cal O_{\Bbb P^{1}}(-1)
\endalign
$$
are isomorphisms. The sheaf $\Cal O_{\Bbb P^{1}}(-1)$ has a metric
induced by the standard metric on $\Bbb C^{2}$. We put on $\Cal
I_{x_{i}}$ and $\Cal I_{y_{i}}$ the metrics induced by the above
isomorphisms. By 2.9, these metrics are smooth at infinity.

For each pair
of integers $i\in \{1,\dots ,n\}$ and $j\in \{-1,0\}$, we write
$$
\Cal I_{i,j}=\cases
\Cal I_{y_{i}},&\text{ if }j=-1,\\
\Cal I_{x_{i}},&\text{ if }j=0.\\
\endcases
$$
For each $\alpha \in \left<-1,0,1\right>^{n}$, with $\alpha \le 0$, and
for each $k\in \{1,\dots ,n\}$, with $\alpha _{k}=-1$, we write
$$
\align
\Cal J_{\alpha }&=\prod_{i=1}^{n} \Cal I^{-1}_{i,\alpha _{i}}\subset
\Cal K_{X\times (\Bbb P^{1})^{n}},\\
\Cal J_{\alpha,k }&=\prod_{i\not =k} \Cal I^{-1}_{i,\alpha
_{i}}\subset \Cal K_{X\times (\Bbb P^{1})^{n}},
\endalign
$$
where $\Cal K_{X\times (\Bbb P^{1})^{n}}$ is the sheaf of rational
functions on $X\times (\Bbb P^{1})^{n}$.

Given an $n$-tuple $\alpha \le 0$ and an integer $k\in \{1,\dots ,n\}$,
with $\alpha _{k}=-1$, we write $\alpha +k=(\alpha _{1},\dots
,\alpha _{k}+1,\dots ,\alpha _{n} )$. We have the inclusions
$$
\align
\Cal J_{\alpha ,k}&\subset \Cal J_{\alpha },\\
\Cal J_{\alpha ,k}&\subset \Cal J_{\alpha +k}.
\endalign
$$
Let us denote by $\varphi _{\alpha ,k}:\overline {\Cal F}_{\alpha
}\longrightarrow \overline {\Cal F}
_{\alpha +k}$ the morphism $\overline {\Cal F}(\alpha \longrightarrow
\alpha +k)$. Let $\psi $ be the morphism
$$
\psi :\bigoplus_{\alpha \le 0}\,\bigoplus_{k|\alpha
_{k}=-1}p_{0}^{*}\overline {\Cal F}
_{\alpha }\otimes \Cal J_{\alpha ,k} \longrightarrow
\bigoplus _{\alpha \le 0}p_{0}^{*}\overline {\Cal F}_{\alpha }\otimes
\Cal J_{\alpha },
$$
which sends $s\otimes g\in p_{0}^{*}\overline {\Cal F}_{\alpha
}\otimes \Cal J_{\alpha ,k}$ to
$$
\align
\psi (s\otimes g)&=s\otimes g+\varphi _{\alpha ,k}(s)\otimes g\\
&\in
p_{0}^{*}\overline {\Cal F}_{\alpha }\otimes \Cal J_{\alpha }\oplus
p_{0}^{*}\overline {\Cal F}_{\alpha+k }\otimes \Cal J_{\alpha+k }.
\endalign
$$

The locally free sheaf
$
\bigoplus _{\alpha \le 0}p_{0}^{*}\overline {\Cal F}_{\alpha }\otimes
\Cal J_{\alpha }$
has a metric induced by the metrics of $\Cal I_{x_{i}}$, $\Cal
I_{y_{i}}$ and $\overline {\Cal F}_{\alpha }$. This metric is smooth
at infinity.

\sec
{\bf Definition 3.8.} The $n$-transgression of $\overline {\Cal F}$
is the hermitian locally free
sheaf
$$
\tr_{n}(\overline {\Cal F})=\Coker (\psi ),
$$
with the metric induced by the metric of
$
\bigoplus _{\alpha \le 0}p_{0}^{*}\Cal F_{\alpha }\otimes \Cal J_{\alpha
}$.
By proposition 2.8, this metric is smooth at infinity.

The following result follows directly from the definition.

\proclaim{Proposition 3.9} Let $\overline {\Cal F}$ be an emi-$n$-cube.
Then there are isometries
$$
\align
\tr_{n}(\overline {\Cal F})|_{\{x_{i}=0\}}&\cong
tr_{n-1}(\partial ^{0}_{i}\overline {\Cal F}),\\
\tr_{n}(\overline {\Cal F})|_{\{y_{i}=0\}}&\cong
tr_{n-1}(\partial ^{-1}_{i}\overline {\Cal F})\ort
tr_{n-1}(\partial ^{1}_{i}\overline {\Cal F}).
\endalign
$$
\endproclaim

\sec (3.7) Let us give an inductive construction of the transgressions.
If $n=1$, an emi-$1$-cube, $\overline {\Cal F}$ is a
short exact sequence
$$
\overline {\Cal F}_{-1}@>f>>
\overline {\Cal F}_{0}@>>>
\overline {\Cal F}_{1},
$$
where the metric of $\overline {\Cal F}_{1}$ is induced by the metric
of $\overline {\Cal F}_{0}$. Then $\tr_{1}(\overline {\Cal F})$ is the
cokernel of the map
$$
\matrix
\Cal F_{-1}&\longrightarrow&
\Cal F_{-1}\otimes \Cal I_{y_{1}}^{-1}&\oplus&
\Cal F_{0}\otimes \Cal I_{x_{1}}^{-1}\\
s&\longmapsto&s\otimes 1&\oplus& f(s)\otimes 1.
\endmatrix
$$

Observe that this is a minor modification of the locally free sheaf used
by Bismut, Gillet and Soul\'e (\cite {B-G-S}, \cite {G-S 1}) to
construct
Bott-Chern forms. In the definition given here, we avoid the use of
partitions of unity,
obtaining a natural construction. The price is to restrict ourselves to
emi-$n$-cubes.

If $\overline {\Cal F}$ is an emi-$n$-cube, let
$\tr_{1}(\overline {\Cal F})$ be the emi-$n-1$-cube over $X\times
\Bbb P^{1}$ defined by:
$$
\tr_{1}(\overline {\Cal F})_{\alpha }=
\tr_{1}(\partial ^{\alpha }_{n^{c}}\overline {\Cal F}).
$$
Then we write
$$
\tr_{k}(\overline {\Cal F})=\tr_{1}(\tr_{k-1}(\overline {\Cal F})).
$$
The hermitian locally free sheaf $\tr_{n}(\overline {\Cal F})$ defined
in this way coincides with the earlier definition.
Thus the transgressions are simply an iteration of the
construction of Bismut, Gillet and Soul\'e.

\sec
(3.8)
For any homology complex $A_{*}$, we shall denote by $A^{*}$ the
cohomology complex defined by $A^{k}=A_{-k}$.
Let us use the transgressions previously defined, to associate to every
emi-$n$-cube a family of differential forms.

\sec
{\bf Definition 3.10.} Let
$$
\ch:\Bbb ZC_{emi}^{*}(X)\longrightarrow
\bigoplus _{p}\widetilde {\frak H}^{*}(X,p)[2p]
$$
be the map given by
$$
\ch(\overline {\Cal F})=\chw{0}(\tr_{n}(\overline {\Cal
F}))_{\Cal H},
$$
where $\chw{0}(\cdot)_{\Cal H}$ is as in (2.4).

\proclaim{Proposition 3.11} The map $\ch$ is a morphism of complexes and
factorizes through a unique morphism
$$
\ch:\widetilde {\Bbb Z}C_{emi}^{*}(X)\longrightarrow
\bigoplus _{p}\widetilde {\frak H}^{*}(X,p)[2p].
$$
\endproclaim
\demo{Proof} To see that it is a morphism of complexes, observe that,
since the forms $\chw{0}(\cdot)_{\Cal H}$ are closed,
$$
\align
d \ch(\tr_{n}(\overline {\Cal F}))&=
\sum_{i=1}^{n}\sum_{j=0}^{1}
(-1)^{i+j}d^{j}_{i}\chw{0}(\tr_{n}(\overline {\Cal F}))_{\Cal
H}\\
&=\sum_{i=1}^{n} (-1)^{i}\chw{0}(\tr_{n}(\overline {\Cal F}))_{\Cal
H}|_{\{x_{i}=0\}}+
\sum_{i=1}^{n} (-1)^{i+1}\chw{0}(\tr_{n}(\overline {\Cal F}))_{\Cal
H}|_{\{y_{i}=0\}}\\
&=\sum_{i=1}^{n} (-1)^{i}\chw{0}(\tr_{n}(\overline {\Cal F})
|_{\{x_{i}=0\}})_{\Cal H}+
\sum_{i=1}^{n} (-1)^{i+1}\chw{0}(\tr_{n}(\overline {\Cal F})
|_{\{y_{i}=0\}})_{\Cal H}.
\endalign
$$
Therefore, by proposition 3.9,
$$
\align
d \ch(\tr_{n}(\overline {\Cal F}))&=
\sum_{i=1}^{n}\sum_{j=-1}^{1}
(-1)^{i+j}\chw{0}(\tr_{n-1}(\partial ^{j}_{i}\overline {\Cal F}))_{\Cal
H}\\
&=\ch(d\overline {\Cal F}).
\endalign
$$

To see the existence of the factorization,
we have to show that, for a degenerate emi-$n$-cube
$\overline {\Cal F}$,  we have $\ch(\overline {\Cal F})=0$ in
$\bigoplus \widetilde {\frak H}^{*}(X,p)$. By symmetry we may assume
that $\overline {\Cal F}=s^{j}_{n}\overline {\Cal G}$, with $j\in
\{-1,1\}$ and $\overline {\Cal G}$ an emi-$n-1$-cube.

If $j=1$, then $\tr_{n-1}(\overline {\Cal F})$ is the exact sequence
$$
0@>>>
(\Id\times s^{n})^{*}\tr_{n-1}(\overline {\Cal G})
@>\Id >>
(\Id\times s^{n})^{*}\tr_{n-1}(\overline {\Cal G})
@>>>0
@>>>0.
$$
Therefore $\tr_{n}(\overline {\Cal F})$ is the cokernel of the map
$$
\matrix
(\Id\times s^{n})^{*}\tr_{n-1}(\overline {\Cal G})
&\longrightarrow&
(\Id\times s^{n})^{*}\tr_{n-1}(\overline {\Cal G})
\otimes \Cal I_{y_{n}}^{-1}&\oplus&
(\Id\times s^{n})^{*}\tr_{n-1}(\overline {\Cal G})
\otimes \Cal I_{x_{n}}^{-1}\\
x&\longmapsto&x\otimes 1&\oplus& x\otimes 1.
\endmatrix
$$
But $\Cal I_{y_{n}}^{-1}$ and $\Cal I_{x_{n}}^{-1}$ are both isometric
with $p_{n}^{*}\Cal O(1)$. Hence this cokernel is isometric with
$(\Id\times s^{n})^{*}\tr_{n-1}(\overline {\Cal G})\otimes p_{n}^{*}\Cal
O(2)$, where $\Cal O(2)$ is provided with the standard metric. Thus
$$
\chw{0}(\tr_{n}(\overline {\Cal F}))_{\Cal H}=
(\Id\times s^{n})^{*}\chw{0}(\tr_{n-1}(\overline {\Cal G}))_{\Cal H}+
2\omega _{n}\land
(\Id\times s^{n})^{*}\chw{0}(\tr_{n-1}(\overline {\Cal G}))_{\Cal H}.
$$
which is zero in $\bigoplus \widetilde {\frak H}^{*}(X,p)$.

The case $j=-1$ is analogous.
\enddemo

\sec
{\bf Definition 3.12.} We shall denote also by $\ch$ the composition
$$
\widetilde {\Bbb Z}C^{*}(X)@>\lambda  >>
\widetilde {\Bbb Z}C^{*}_{emi}(X)@>\ch >>
\bigoplus _{p}\widetilde {\frak H}^{*}(X,p)[2p].
$$

\sec
{\bf Definition 3.13.} Let $\overline {\Cal F}$ be a metrized exact
$n$-cube. The form $\ch(\lambda (\overline {\Cal F}))$ will be called
the Bott-Chern form of $\overline {\Cal F}$ and will be denoted by
$\chw{n}(\overline {\Cal F})_{\Cal H}$.

\heading {\S}4 Higher characteristic classes.
\endheading

The Chern character from $K$-theory to a suitable cohomology theory,
such as absolute Hodge cohomology, is additive for exact sequences.
Nevertheless, given a cochain complex which computes absolute Hodge
cohomology, we cannot make a consistent choice of representatives of
the Chern character that behaves additively.
Following the ideas of Schechtman (\cite {Sch}), the lack
of additivity at the level of complexes, of the Chern character for
$K_{n}$, gives us the Chern character for $K_{n+1}$.

In the previous section we have associated, to each metrized exact
$n$-cube, a family of differential forms. The differential form
associated to an $n$-cube measures the lack of additivity of the
differential forms associated to its faces.
In this section we shall see that this construction allows us to define
higher Chern character
classes from $K$-theory to absolute Hodge cohomology.

\sec (4.1)
Let us begin by reviewing the Waldhausen $K$-theory of a small exact
category.
We shall follow \cite {Sch} (See also \cite {Wal} or \cite {Lo 1}).

For $n\in \Bbb N$, let $\Cat (n)$ denote the category associated with
the ordered set $\{1,\dots ,n\}$. Let $M_{n}$ be the category of
morphisms of $\Cat(n)$. That is
$$
\Ob M_{n}=\left\{(i,j)\in \Bbb N\times \Bbb N\mid 0\le i\le j\le n
\right\}, $$
and $\Hom ((i,j),(k,l))$ contains a unique element if $i\le k$ and $j
\le l$ and is empty otherwise. The categories $M_{n}$ form a
cosimplicial category $M$.

For any category $\frak C$, let us denote by $M_{n}\frak C$ the category
of functors from $M_{n}$ to $\frak C$.

\sec
{\bf Definition 4.1.} Let $\frak E$ be a small exact category and $0$
a fixed zero object of $\frak E$.
Let $\underline S_{n}\frak E$ be the full subcategory of $M_{n}\frak E$,
whose objects are the functors $M_{n}\longrightarrow \frak E$, such
that,
\roster
\item for all $i$, $E_{i,i}=0$;
\item for all $i\le j\le k$,
$$
E_{i,j}@>>> E_{i,k} @>>> E_{j,k}
$$
is a short exact sequence.
\endroster
Let us write $S_{n}\frak E=\Ob \underline S_{n}\frak E$.
We shall denote by $\underline S\frak E$ or $\underline S^{1}\frak E$
the simplicial
exact category $\coprod \underline S_{n}\frak E$, and by $S\frak E$ or
$S^{1}\frak E$ the simplicial set $\Ob \underline S\frak E$.

In other words we have:
$$
\align
S_{0}\frak E&=\{0\},\\
S_{1}\frak E&=\Ob \frak E,\\
S_{2}\frak E&=\{\text{exact sequences of $\frak E$} \},\\
S_{n}\frak E&=\left \{
\matrix
\text {sequences of monomorphisms}\\
E_{0,1}\rightarrow E_{0,2}\rightarrow\dots \rightarrow E_{0,n}\\
\text {with a choice of quotients}\\
E_{i,j}\cong E_{0,j}/E_{0,i}
\endmatrix
\right\}.
\endalign
$$
In particular $S\frak E$ is a pointed simplicial set. In the sequel we
shall sometimes use the word space to denote simplicial sets.

For a space $C$, we shall denote by $|C|$ its geometric
realization.

\proclaim{Proposition 4.2} {(\rm Cf. \cite {Lo 2}.)} There is a homotopy
equivalence $$
S\frak E \cong BQ\frak E,
$$
where $Q$ denotes Quillen's $Q$-construction and $B$ means classifying
space. Therefore, for all
$i\ge 0$, $$
K_{i}(\frak E)=\pi _{i+1}\left(\left |S_{\cdot}\frak E\right |,0\right).
$$
\endproclaim

\sec
(4.2)
Let us recall the notion of spectrum from \cite {Th}.
For any pointed space $C$, let us write
$\Sigma C$ for the suspension of $C$, and $\Omega
C$ for the loop space of $C$. We shall use
the same notation for topological spaces.

\sec
{\bf Definition 4.3.} A prespectrum $X$ is a sequence of pointed
spaces $X_{n}$ for all non-negative integers $n$, together with
structure maps $\Sigma X_{n}\longrightarrow X_{n+1}$. These maps can
also
be described by their adjoint $X_{n}\longrightarrow \Omega X_{n+1}$. A
fibrant spectrum is a prespectrum such that all $X_{n}$ are fibrant
spaces and the structure maps $X_{n}\longrightarrow \Omega X_{n+1}$ are
weak equivalences.

The space $S\frak E$ is a piece of a
prespectrum. To construct the other spaces that
form the prespectrum,
we write inductively
$$
\align
\underline S^{m}\frak E&= \underline S \underline
S^{m-1} \frak E, \\
S^{m}\frak E&= S\underline
S^{m-1}\frak E. \\
\endalign
$$
Then $\underline S^{m}$ is an exact $m$-simplicial category and
$S^{m}$ is a $m$-simplicial set. For a poly-simplicial set $C$ let
$\diag (C)$ denote its diagonal space. We shall denote by
$\left|C\right|=\left|\diag(C)\right|$ its geometric realization.

\proclaim{Proposition 4.4} {\rm (Schechtman  \cite {Sch 1.2}.)} There
are natural maps $$
\varphi _{m}:\Sigma S^{m}\frak E \longrightarrow
S^{m+1}\frak E,
$$
inducing homotopy equivalences
$$
\left| S^{m}\frak E\right |
\cong \Omega \left |S^{m+1}\frak E\right |.
$$
\endproclaim

As a consequence of this proposition, if we write $S^{0}\frak E=\Omega
S\frak E$, then the sequence of spaces $\diag(S^{m}\frak E)$
is a prespectrum. Moreover, if we replace the above spaces by
weakly equivalent fibrant spaces we shall obtain a fibrant spectrum. For
instance, let us denote by $\Sing$ the singular functor (see \cite
{B-K}).
Then, if we write
$$
\bold{K}_{m}(\frak E)=\Sing(|S^{m}\frak E|),
$$
the spaces $\bold{K}_{m}$ form a fibrant spectrum. By proposition
4.2, the homotopy of this fibrant spectrum is the $K$-theory of $\frak
E$.

\sec
(4.3) For example,
let $X$ be a smooth complex variety, let $m \ge 1$ be an integer and let
us write $S^{m}(X)=S^{m}(\overline {\frak E}(X))$.
Then, for  $i\ge 0$, the $i$-th $K$-group of $X$ is
$$
K_{i}(X)=\pi _{i+m}(|S^{m}(X)|,0).
$$
By proposition 4.4 this definition does not depend on the choice of $m$.

\sec
(4.4)
Let us associate, to each element of $S_{n}\frak E$, an exact
$n-1$-cube. We shall do so inductively. For $n=1$, we write
$$
\Cub(\{E_{i,j}\}_{0\le i\le j\le 1})= E_{0,1}.
$$
Assume that we have defined $\Cub E$ for all $E\in S_{m}\frak E$,
with $m<n$. Let $E\in S_{n}\frak E$. Then $\Cub E$  is
the $n-1$-cube with
$$
\align
\partial ^{-1}_{1}\Cub E&=s^{1}_{n-2}\dots s^{1}_{1}(E_{0,1}),\\
\partial ^{0}_{1}\Cub E&=\Cub (\partial _{1}E),\\
\partial ^{-1}_{1}\Cub E&=\Cub(\partial _{0}E).\\
\endalign
$$
For instance, if $n=2$, then $\Cub \left(\{E_{i,j}\}_{0\le i\le j\le 2
}\right)$ is the short exact sequence
$$
E_{0,1}@>>> E_{0,2}@>>> E_{1,2}.
$$
On the other hand, if $n=3$, then
$\Cub \left(\{E_{i,j}\}_{0\le i\le j\le 3
}\right)$ is the exact square
$$
\CD
E_{0,1}@>>> E_{0,2}@>>> E_{1,2}\\
@VVV @VVV @VVV\\
E_{0,1}@>>> E_{0,3}@>>> E_{1,3}\\
@VVV @VVV @VVV\\
0@>>> E_{2,3}@>>> E_{2,3}\\
\endCD
$$

All the faces of the $n-1$-cube $\Cub E$ can be computed explicitly.

\proclaim{Proposition 4.5} Let $E\in S_{n}\frak E$. Then, for
$i=1,\dots ,n-1$, the faces of the $n-1$-cube $\Cub E$ are
$$
\align
\partial ^{-1}_{i}\Cub E&= s^{1}_{n-2}\dots s^{1}_{i}\Cub
\partial _{i+1}\dots \partial _{n}E,\\
\partial ^{0}_{i}\Cub E&= \Cub
\partial _{i}E,\\
\partial ^{1}_{i}\Cub E&= s^{-1}_{i-1}\dots s^{-1}_{1}\Cub
\partial _{0}\dots \partial _{i-1}E.\\
\endalign
$$
\endproclaim

By proposition 4.5 and using induction we have,

\proclaim{Corollary 4.6} The $n-1$-cube $\Cub E$ is exact.
\endproclaim

Therefore we have a map $\Cub:S_{n}\frak E\longrightarrow C_{n-1}\frak E$.

\sec
(4.5)
Let $\Bbb ZS\frak E$ be the homological complex associated with the
simplicial set $S\frak E$. That is, $\Bbb ZS_{n}\frak E$ is the free
abelian group generated by $S_{n}\frak E$, and the differential
$d:\Bbb ZS_{n}\frak E\longrightarrow \Bbb ZS_{n-1}\frak E$ is given by
$$
d=\sum_{i=0}^{n}(-1)^{i}\partial _{i}.
$$

The map $\Cub$ can be extended by linearity to a map
$$
\Cub:\Bbb ZS\frak E[1]\longrightarrow \Bbb ZC\frak E.
$$
Note that this map is not a morphism of complexes.
However, the map $\Cub$ induces a map also denoted by
$\Cub:\Bbb ZS\frak E[1]\longrightarrow \widetilde {\Bbb Z}C\frak E$.
And, since by proposition  4.5,
$$
d\Cub E=\Cub dE+\text{degenerate elements},
$$
we have:

\proclaim{Corollary 4.8} The map
$\Cub:\Bbb ZS\frak
E[1]\longrightarrow \widetilde {\Bbb Z}C\frak E$ is a morphism of
complexes.
\endproclaim

\sec
(4.6)
We can obtain analogous maps for all the spaces $S^{m}\frak E$. In
particular, we have maps
$$
\Cub: S_{n_{1}}\dots S_{n_{m}}\frak E\longrightarrow
C_{n_{1}-1}\dots C_{n_{m}-1}\frak E\longrightarrow
C_{n_{1}+\dots +n_{m}-m}\frak E.
$$
Let us denote by $\Bbb ZS^{m}\frak E$ the chain complex that, in degree
$n$, is the free abelian group generated by
$$
\coprod_{n_{1}+\dots +n_{m}=n}S_{n_{1}}\dots S_{n_{m}}\frak E.
$$
The differential of this complex is the alternate sum of all the face
maps. Note that this complex is homotopically equivalent to $\Bbb
Z\diag(S^{m}\frak E)$. The induced map
$$
\Cub: \Bbb ZS^{m}\frak E[m]\longrightarrow
\widetilde {\Bbb Z}C \frak E
$$
is also a morphism of complexes.

\sec
(4.7)Let $m\ge 1$ be an integer.
We shall denote by $\Bbb ZS^{*}_{m}(X)$ the cohomological complex
associated to the homological complex $\Bbb ZS_{*}^{m}(X)$.

\sec
{\bf Definition 4.9.} The Chern character map is the composition
$$
\Bbb ZS^{*}_{m}(X)[-m]@>\Cub >>
\widetilde {\Bbb Z}C^{*} (X)@>\lambda >>
\widetilde {\Bbb Z}C^{*}_{emi}(X)@>\ch >>
\bigoplus _{p}\widetilde {\frak H}^{*}(X,p)[2p].
$$
This map will also be denoted by $\ch$. The Chern character classes are
obtained by composing with the Hurewicz map:
$$
K_{i}(X)=\pi _{i+m}(S^{m}(X)))@>>>
H_{i+m}(\Bbb ZS^{m}(X))@>>>
\bigoplus_{p}H^{2p-i}_{\Cal H}(X,p).
$$

\proclaim{Proposition 4.10} The above definition does not depend on the
choice of $m$. \endproclaim
\demo{Proof} For a pointed simplicial set $S$, with base point $p$,
we shall write
$$
\Bbb Z'S_{*} = \Bbb Z S_{*}/\Bbb Z p_{*}.
$$
The natural map $\Sigma
S^{m}(X)\longrightarrow S^{m+1}(X)$ induces a morphism
$$
\Bbb Z' S^{m}_{*}(X) [m] \longrightarrow
\Bbb Z' S^{m+1}_{*}(X) [m+1].
$$
By the proof of proposition 4.4 in \cite {Sch 1.2} this map is induced
by
the natural bijection $S^{m}(X)\longrightarrow S_{1}S^{m}(X)$. Therefore
the diagram
$$
\CD
\Bbb Z' \Sigma S^{m}_{*}(X) [m+1] @>>> \Bbb Z' S^{m+1}_{*}(X)[m+1]\\
@AAA @VV\Cub V\\
\Bbb Z' S^{m}_{*}(X) [m] @>\Cub >> \widetilde {\Bbb Z}'C_{*} (X)
\endCD
$$
is commutative. By the commutativity of this diagram and proposition 4.4
the Chern character classes are independent of the choice of $m$.

\enddemo

\heading {\S}5 Beilinson's regulator.
\endheading

The aim of this section is to prove that the higher Chern character
classes defined in {\S}4 agree with Beilinson's
regulator map.

\sec
(5.1) Let us begin by extending the definition of the map $\ch$ to the
case of simplicial smooth complex varieties. To this end, we first
recall the construction of the absolute
Hodge cohomology of $X=X_{\cdot}$ a smooth simplicial complex variety.
For each $p$, the complexes $\widetilde {\frak H}^{*}(X_{n},p)$ form a
cosimplicial complex as $n$ varies. Let
$\Cal N \widetilde {\frak H}^{*}(X_{\cdot},p)$ be the associated double
complex and let us denote the simple complex by
$$\widetilde {\frak H}^{*}(X,p)=
s(\Cal N \widetilde {\frak H}^{*}(X_{\cdot},p)).
$$
Then
$$
H^{*}_{\Cal H}(X_{\cdot},\Bbb R(p))=
H^{*}(\widetilde {\frak H}^{*}(X,p)).
$$

For the definition of $K$-theory of
simplicial schemes we shall follow \cite {Sch}.
We shall say that a smooth simplicial scheme $X=X_{\cdot}$ has finite
dimension if there is an integer $m$ such that
$$
X=\Sk_{m}(X),
$$
where $\Sk_{m}(X)$ is the $m$-th skeleton of $X$, that is, the
simplicial scheme generated by $X_{0},\dots ,X_{m}$.

Let $X=X_{\cdot}$ be a simplicial
scheme of finite dimension. The family of prespectrums
$\{S(X_{n})\}_{n}$ form a cosimplicial
prespectrum $S(X_{\cdot})$. Let $\bold K(X_{\cdot})$ be a fibrant
cosimplicial fibrant spectrum weakly equivalent to $S(X_{\cdot})$.
Then the $K$ groups of $X$ are defined as
$$
K_{i}(X)=\pi _{i+1}(\Tot \bold K(X_{\cdot})).
$$
Since $X$ has finite dimension, there is a convergent spectral sequence
$$
E_{1}^{p,q}=K_{-q}(X_{p})\Longrightarrow K_{-p-q}(X_{\cdot}).
$$

Observe that for a given simplicial scheme $X_{\cdot}$, of finite
dimension, it is not necessary to work with the whole spectrum.
Let $m$ be such that $X=\Sk_{m}X$.
Let us choose a positive integer $q$, and let $\bold{K}_{q}(X_{\cdot})$
be a fibrant cosimplicial fibrant space, weakly equivalent to
$S^{q}(X_{\cdot})$. If $q>m$ or $q>-i$, then
$$
K_{i}(X_{\cdot})=\pi _{i+q}(\Tot \bold{K}_{q}(X_{\cdot})).
$$

For an arbitrary simplicial scheme
we write
$$
\widehat K_{*}(X_{\cdot})=\lim\Sb \longleftarrow\\m\endSb
K_{*}(\Sk_{m}(X_{\cdot})).
$$

Let $X$ be a smooth simplicial complex variety of finite dimension.
Since the map $\ch$ defined in section {\S}4 gives us a morphism of
complexes
$$
\ch:s\Cal N \Bbb Z S^{*}_{q}(X_{\cdot})[q]\longrightarrow
\bigoplus _{p}\widetilde {\frak H}^{*}(X,p)[2p],
$$
we can extend the definition of the Chern character to
the simplicial case, obtaining maps:
$$
\ch:K_{i}(X)\longrightarrow \bigoplus_{p}H^{2p-i}_{\Cal H}(X,\Bbb R(p)).
$$

If $X$ does not have finite dimension, taking limits, we have also
characteristic classes
$$
\ch:\widehat K_{i}(X)\longrightarrow \bigoplus_{p}H^{2p-i}_{\Cal
H}(X,\Bbb R(p)).
$$

\sec
{\bf Remark. 5.1.} All the constructions needed to define the map
$\ch$ can be extended to the case of a smooth
simplicial scheme $X_{\cdot}$ over $\Bbb
C$, such that each $X_{n}$ is a (not necessarily finite)
disjoint union of
smooth complex varieties. For instance, by a compactification of $X_{n}$
we shall mean a disjoint union of compactifications of each component of
$X_{n}$.

\sec
(5.2)
Beilinson (\cite {Be}) has defined characteristic classes from
$K$-theory to absolute Hodge cohomology. These classes are a
particular case of the characteristic classes defined by Gillet (\cite
{Gi}) to any suitable cohomology theory. In particular, Beilinson's
regulator is the Chern character in this theory. Let us denote
by $\rho $ the Beilinson's regulator map.

Then $\ch$ and $\rho $ are natural
transformations between
contravariant functors. Both agree with the classical Chern character
on the $K_{0}$ groups of smooth complex varieties.
The aim of this section is to prove the following theorem.

\proclaim{Theorem 5.2} Let $X$ be a smooth complex variety.
Let $\sigma \in K_{i}(X)$. Then $\ch(\sigma )=\rho (\sigma )$.
\endproclaim
\demo{Proof}
Let
$\frak U=\{U_{\alpha }\}$ be an open covering of $X$. We shall denote by
$\frak E(X,\frak U)$ the full subcategory of $
\frak E(X)$ composed by the locally free sheaves on $X$
whose restrictions to all $U_{\alpha }$ are free. We shall denote by
$\overline {\frak E}(X,\frak U)$ the category of hermitian vector
bundles on $X$ whose restrictions to all $U_{\alpha }$ are free. Let us
write
$$
K_{i}(X,\frak U)=\pi _{i+1}(S\frak E(X,\frak U))
=\pi _{i+1}(S\overline {\frak E}(X,\frak U)).
$$
Then
$$
K_{i}(X)=\lim_{\Sb\longrightarrow\\\frak U\endSb } K_{i}(X,\frak U).
$$

\sec
(5.2.1)
Following Schechtman (\cite {Sch}) we know that there is a simplicial
scheme $\Cal B P$, which is a classifying space for algebraic
$K$-theory. More precisely, Schechtman proves the following result.

\proclaim{Theorem 5.3} {\rm (Schechtman)} There is a homotopy
equivalence
$$
S\frak E(X,\frak U)\cong \underline {\Hom} (N\frak U,\Cal BP),
$$
where $\underline {\Hom}$ is the function space and $N\frak U$ is the
nerve of the covering.
\endproclaim

\sec
(5.2.2)
Let $Y=Y_{\cdot}$ be a smooth simplicial scheme of finite dimension. Let
us denote by $\underline {\underline{\Hom}}(Y,\Cal BP)$ the cosimplicial
simplicial set $\underline {\underline{\Hom}}(Y,\Cal BP)^{n}_{m}=\Hom
(Y_{n},\Cal BP_{m})$. Then
$$
\underline {\Hom}(Y,\Cal BP)=\Tot \underline{\underline{\Hom}}(Y,\Cal
BP). $$
For any scheme $X$, let $\sigma X$ be the simplicial scheme with
$\sigma X_{n}=X$ and all the faces and degeneracies equal to the
identity. Then the simplicial set $\Hom (X,\Cal BP)$ is the function
space $\underline {\Hom}(\sigma X,\Cal BP)$. Observe that $\sigma X$ is
the nerve of the trivial covering $\{X\}$. Thus, by Theorem 5.3 and
the comparison between $\Tot$ and $ \Holim$ (cf. \cite {B-K, XI, 4.4},
\cite {Th, 5.25} and \cite {Le, 3.1.2}) we
obtain a natural map 
$$
\Tot \underline {\underline
{Hom}}(Y,\Cal BP) \longrightarrow \Holim  S\frak E(Y_{n},\{Y_{n}\}).
$$
Taking homotopy groups we obtain natural maps
$$
\pi _{i} \underline {\Hom}(Y,\Cal BP) \longrightarrow
K_{i-1}(Y).
$$

In particular, if $f:Y \longrightarrow
\Cal BP $ is a simplicial morphism of simplicial schemes, then $f$
defines an element
of $\pi _{0}\underline {\Hom}(Y,\Cal BP)$. Let us denote by $e_{f}$ the
image of this element in $K_{-1}(Y)$.

If $Y$ does not have finite dimension, since any simplicial morphism
$f:Y\longrightarrow \Cal BP$ induces simplicial morphisms
$f_{m}:\Sk_{m}Y\longrightarrow \Cal BP$, taking a limit, we obtain an
element $e_{f}\in \widehat K_{-1}(Y)$.

\sec
{\bf Remark 5.4.} The identity of $\Cal BP$ defines an element,
denoted by $e_{_{\Cal BP}}\in \widehat K_{-1}(\Cal BP)$. Moreover, for
$f$ as above, $e_{f}=f^{*}(e_{_{\Cal BP}})$.

\sec
(5.2.3) The element $e_{_{\Cal BP}}$ is, in some sense, a universal
element in $K$-theory. Since $e_{_{\Cal BP}}\in \widehat K_{-1}$, to
exploit the
universality of this element, we need to relate elements in
$\widehat K_{n}$ with elements in $\widehat K_{-1}$. This can be done
using spheres.

Let $\sigma \in K_{n}(X)$. Then there is
an
open covering $\frak U$ of $X$, such that $\sigma \in K_{n}(X,\frak U)=
\pi _{n+1}(\underline {\Hom}(N\frak U,\Cal BP))$. Therefore, there
exists an integer $d\ge 0$ such that $\sigma $ is
represented by an element
$$
\gamma _{\sigma }\in
\underline {\Hom}(\Sub_d S^{n+1}\times N\frak U,\Cal BP)
=\underline {\Hom}(\Sub_d S^{n+1},
\underline {\Hom}(N\frak U,\Cal BP)),
$$
where $S^{n+1}$ is the (pointed) simplicial $n+1$-dimensional sphere
and $\Sub_d$ is the $d$-th subdivision. Let us denote by
$\Sigma^{n+1}=\Sub_d S^{n+1}$

\proclaim{Lemma 5.5} Let $Y=Y_{\cdot}$ be a smooth simplicial complex
variety. Then there are natural decompositions
$$
\align
\widehat K_{-1}(\Sigma ^{n+1}\times Y)&=\widehat K_{-1}(Y)\oplus
\widehat K_{n}(Y),\\
H^{2p+1}_{\Cal H}(\Sigma ^{n+1}\times Y,\Bbb R(p))&=
H^{2p+1}_{\Cal H}(Y,\Bbb R(p))\oplus
H^{2p-n}_{\Cal H}(Y,\Bbb R(p)).
\endalign
$$
Moreover, the maps $\ch$ and $\rho $ are compatible with these
decompositions. \endproclaim
\demo{Proof} We may assume that $Y$ has finite dimension because the
general case is obtained taking a limit. Then
$$
K _{-1}(\Sigma ^{n+1}\times Y)=\pi _{0}(\Tot_{\alpha }(\Tot
_{\beta }(\bold{K}(\Sigma ^{n+1}_{\alpha }\times Y_{\beta })))).
$$
The spectral sequence associated with $\Tot _{\alpha }$ has
$E_{2}$-term:
$$
E_{2}^{p,q}=\cases K_{-q}(Y),&\text{ if } p=0,n+1,\\
                      0,&\text{ if } p\not=0,n+1.
\endcases
$$
Let us denote by $*$ the simplicial point. Since the spectral sequence
of $*\times Y$ splits the spectral sequence of $\Sigma ^{n+1}\times Y$, the
above spectral sequence degenerates at the $E_{2}$-term, and the exact
sequence obtained from this spectral sequence splits in a natural way.

The same argument works for cohomology. Moreover, since $\ch$ and $\rho $
are natural transformations, they induce morphisms between the
$K$-theoretical and the cohomological spectral sequences, proving the
compatibility statement. \enddemo

Let us denote by $pr:K_{-1}(\Sigma ^{n+1}\times N\frak U)\longrightarrow
K_{n}(N \frak U)$ the projection. The precise meaning of the
universality of $e_{_{\Cal BP}}$ is given by the following result.

\proclaim{Lemma 5.6} In the group $K_{n}(N\frak U)$, the
equality $$
\pr(\gamma _{\sigma }^{*}(e_{_{\Cal BP}}))=\sigma
$$
holds.
\endproclaim
\demo{Proof} By remark 5.4,
$$
\pr(\gamma _{\sigma }^{*}(e_{_{\Cal BP}}))=\pr(e_{\gamma _{\sigma }}).
$$
On the other hand, by the definition of $\gamma _{\sigma }$, the map
$$
\pi _{0}(\underline {\Hom}(\Sigma ^{n+1}\times N\frak U,\Cal BP))
@>\pr >>
\pi _{n+1}(\underline {\Hom}(N\frak U,\Cal BP))
$$
sends the class of $\gamma _{\sigma }$ to the class of $\sigma $.
Therefore, since the diagram
$$
\CD
\pi _{0}(\underline {\Hom}(\Sigma ^{n+1}\times N\frak U,\Cal BP))
@>>> K_{-1}(\Sigma ^{n+1}\times N\frak U)\\
@V\pr VV @V\pr VV\\
\pi _{n+1}(\underline {\Hom}(N\frak U,\Cal BP))
@>>> K_{n}(N\frak U).
\endCD
$$
is commutative we have that
$$
\pr(\gamma _{\sigma }^{*}(e_{_{\Cal BP}}))=\sigma .
$$
\enddemo

\sec
(5.2.4)
By Remark 5.1, the map $\ch$ is defined for the simplicial
scheme $\Cal BP$. Moreover,
by the naturality of $\ch$ and $\rho $ and their compatibility with the
map $\pr$, we have
$$
\align
\ch(\sigma )&=\pr(\gamma ^{*}_{\sigma }(\ch(e_{_{\Cal BP}}))),\\
\rho (\sigma )&=\pr(\gamma ^{*}_{\sigma }(\rho (e_{_{\Cal BP}}))).
\endalign
$$
Thus, to prove theorem 5.2, we are led to compare $\ch(e_{_{\Cal BP}})$
and $\rho (e_{_{\Cal BP}})$. For this comparison, we need to understand
the cohomology of $\Cal BP$. This cohomology has been computed by
Schechtman (\cite {Sch}).
The simplicial scheme $\Cal BP $ is the classifying space of a
simplicial
group $P_{\cdot}$, where $P_{0}=*$ and $P_{1}=\coprod_{n}GL(n)$. Thus it
is a bisimplicial
scheme $\Cal B_{\cdot}P_{\cdot}$. The edge homomorphism of the spectral
sequence associated to the second index gives us a morphism
$$
d_{H}:H^{2p+1}_{\Cal H}(\Cal BP,\Bbb R(p))\longrightarrow
 \prod_{n\ge 0}  H^{2p}_{\Cal H}(BGL(n),\Bbb R(p)).
$$
Let us denote by $A=H^{*}_{\Cal H}(\Spec \Bbb C,\Bbb R(*))$.

For each $i,n$ let us denote by
$$
c_{i,n}=c_{i}(E_{n})\in H^{2i}_{\Cal H}(BGL(n),\Bbb R(i)),
$$
the $i$-th Chern class of the tautological vector bundle over $BGL(n)$.
Then we have an isomorphism
$$
H^{*}_{\Cal H}(BGL(n),\Bbb R(*))=A[c_{1,n},\dots ,c_{n,n}].
$$
Let $s_{k,n}\in A[c_{1,n},\dots ,c_{n,n}]$ be the $k$-th Newton
polynomial in the $c_{i,n}$. That is, $s_{k,n}/n!$ is the degree $k$
term of the Chern character of the tautological vector bundle $E_{n}$.
Let us write
$$
s_{k}=(s_{k,0},s_{k,1},\dots )\in \prod_{n\ge 0} H^{2k}_{\Cal
H}(BGL(n),\Bbb R(k)).
$$
\proclaim{Proposition 5.7} {\rm (Schechtman \cite {Sch})}
There exist elements $s^{1}_{k}\in H^{2k+1}_{\Cal H}(\Cal BP,\Bbb
R(k))$ such that $d_{H}(s^{1}_{k})=s_{k}$ and
$$
H^{*}_{\Cal H}(\Cal BP,\Bbb R(*))=A[s^{1}_{0},s^{1}_{1},\dots ].
$$
\endproclaim

\sec
(5.2.5)
Since
$$
H^{n}_{\Cal H}(\Spec \Bbb C,\Bbb R(p))\cong
\cases
\Bbb R, & \text{ if } n=p=0, \text{ or } n=1,\  p>0,\\
0,&\text{ otherwise,}
\endcases
$$
any element of $H_{\Cal H}^{2k+1}(\Cal BP,\Bbb R(k))$ can be written as
$$\alpha s^{1}_{k}+ \text{ decomposable elements},
$$
with  $\alpha \in \Bbb R$. Moreover, since by the proof of 5.7 (\cite
{Sch}) the decomposable elements are mapped to $0$ by $d_{H}$, we have

\proclaim{Corollary 5.8}
The group $\Ker d_{H}\subset \bigoplus_{k} H_{\Cal H}^{2k+1}(\Cal
BP,\Bbb R(k))$ is generated by decomposable elements. \endproclaim

\sec
(5.2.6)
Schechtman computes the groups $\widehat K_{*}(\Cal
BP)$ in a similar way. In particular, there is also an edge homomorphism
$$
d_{K}:\widehat K_{-1}(\Cal B P)\longrightarrow
\prod_{n} \widehat K_{0}(\Cal BGL(n)).
$$
Moreover, by the naturality of $\ch$ and $\rho $, they are compatible
with the edge homomorphisms. In particular
$$
d_{H}(\rho (e_{_{\Cal BP}}))=
\rho (d_{K}(e_{_{\Cal BP}})), \text{ and }
d_{H}(\ch (e_{_{\Cal BP}}))=
\ch (d_{K}(e_{_{\Cal BP}})).
$$

\sec
(5.2.7) Our next step will be to compare $d_{H}(\rho (e_{_{\Cal BP}}))$
with $d_{H}(\ch (e_{_{\Cal BP}}))$. To this end we shall
see that, since the maps $\ch$ and $\rho $ agree for the
$K_{0}$ groups
of smooth complex varieties then they also agree for the group
$\widehat K_{0}(BGL(n))$.

\proclaim{Proposition 5.9} Let $\sigma \in \widehat K_{0}(BGL(n))$.
Then
$$
\ch(\sigma )=\rho(\sigma ).
$$
\endproclaim
\demo{Proof} Let $Gr(n,k)$ be the Grassman manifold of dimension $n$
linear subspaces of $C^{k}$ and let $E(n,k)$ be the rank $n$
tautological vector bundle. Let $\frak U_{k}=\{U_{\alpha }\}$  be
the standard trivialization of $E(n,k)$.
Let us denote by $\psi :N\frak
U_{k}\longrightarrow Gr(n,k)$ the natural map and by $\varphi
_{k}:N\frak
U_{k}\longrightarrow BGL(n)$ the classifying map. Since
absolute Hodge
cohomology  can be computed as the cohomology of a Zariski sheaf, the
map
$$
\psi ^{*}:H^{*}_{\Cal H}(Gr(n,k),\Bbb R(*)) \longrightarrow
H^{*}_{\Cal H}(N\frak U_{k},\Bbb R(*))
$$
is an isomorphism.
Moreover, for each $i_{0}$ there is a number $k_{0}$, such that, for all
$k\ge k_{0}$ and all $i\le i_{0}$ the map
$$
\varphi _{k}^{*}:H^{i}_{\Cal H}(BGL(n),\Bbb R(*))\longrightarrow
H^{i}_{\Cal H}(N\frak U_{k},\Bbb R(*))
$$
is an isomorphism. But for $\sigma \in \widehat K_{0}(BGL(n))$ we have
$$
\varphi _{k}^{*}(\ch(\sigma ))=
\ch(\varphi _{k}^{*}(\sigma ))=
\rho(\varphi _{k}^{*}(\sigma ))=
\varphi _{k}^{*}\rho((\sigma )).
$$
Since this is true for all $k$ we have $\ch(\sigma )=\rho (\sigma )$.
\enddemo

Combining 5.8 and 5.9 we get:

\proclaim{Corollary 5.10} The element $\ch(e_{_{\Cal BP}})-\rho
(e_{_{\Cal BP}})$ belongs to $\Ker d_{H}$. Therefore it is
a sum of decomposable elements.\endproclaim

\sec
(5.2.8)
To exploit the fact that
$\ch(e_{_{\Cal BP}})-\rho (e_{_{\Cal BP}})$ is
a sum of decomposable elements, we shall
give a description of how a class in $H^{*}_{\Cal
H}(\Cal BP, \Bbb R(*))$ determines a map between $K$-theory and
absolute Hodge cohomology.

For any smooth simplicial scheme over $\Bbb C$, $X$, and integers $n$,
$p$, the complex
$$
\Cal H^{*}(X,n,p)=\tau _{\le 0}\widetilde {\frak H}^{*}(X,p)[n]
$$
is a negatively graded cohomological complex. Let $\Cal H_{*}(X,n,p)$ be
the associated homological complex. Let us denote
by $\Cal K(X,n,p)$ the simplicial group obtained by Dold-Puppe from
$\Cal H_{*}(X,n,p)$. Then, for $i\ge 0$,
$$
\pi _{i}\Cal K(X,n,p)=H_{\Cal H}^{n-i}(X,\Bbb R(p)).
$$

Let us fix a smooth complex variety $X$, and $\frak U$ an open covering
of $X$. Let us denote by $\varphi $ the tautological map
$$
\varphi :N\frak U\times \underline {\Hom}(N\frak U,\Cal
BP)\longrightarrow \Cal BP.
$$
Given any class $x\in H^{n}_{\Cal H}(\Cal BP,\Bbb R(p))$, we have a
class
$$
\varphi ^{*}(x)\in H^{n}_{\Cal H}(N\frak U\times
\underline {\Hom}(N\frak U,\Cal BP),\Bbb R(p))=
\Hom_{Ho}(\underline {\Hom}(N\frak U,\Cal BP),
\Cal K(N\frak U,n,p)).
$$
For any integer $i$, let us denote by $\pi_{i} (x)$ the induced map
$$
\pi_{i} (x):
K_{i-1}(X,\frak U)=\pi _{i}\underline {\Hom} (N\frak U,\Cal BP)
\longrightarrow \pi _{i}\Cal K(N\frak U,n,p)=H_{\Cal H}^{n-i}(X,\Bbb
R(p)).
$$
Taking the limit over all coverings we obtain morphisms
$$
\pi_{i} (x):
K_{i-1}(X)
\longrightarrow H_{\Cal H}^{n-i}(X,\Bbb
R(p)).
$$
This construction can be extended to the case when $X$ is a simplicial
smooth complex manifold.

\proclaim{Lemma 5.11} For $x\in H^{2k+1}_{\Cal H}(\Cal BP,\Bbb R(p))$
and $\sigma \in K_{i-1}(X,\frak U)$ we have
$$
\pi _{i}(x)(\sigma )=\pr (\gamma_{\sigma } ^{*}(x)),
$$
where $\gamma _{\sigma }$ is as in (5.2.3).
\endproclaim
\demo{Proof} Since the map $\pi _{*}(x)$ is natural, the same argument
as for $\ch$ and $\rho $ shows that
$$
\pi _{i}(x)(\sigma )=\pr (\gamma _{\sigma }^{*}(\pi _{0}(x)(e_{_{\Cal
BP}}))). $$

Let us denote by $\widehat K_{i}(\Cal BP,\Cal BP)$ the $K$-theory groups
of $\Cal BP$ with respect to the trivial covering. Then the map
$$
\pi_{0} (x):
\widehat K_{-1}(\Cal BP,\Cal BP)=\pi _{0}\underline {\Hom} (\Cal BP,\Cal
BP) \longrightarrow \pi _{0}\Cal K(\Cal BP,2k+1,k)=H^{2k+1}(\Cal BP,\Bbb
R(k)),
$$
sends the class of $f\in \pi _{0}\underline {\Hom} (\Cal BP,\Cal BP)$ to
$f^{*}(x)$. Since $e_{_{\Cal BP}}$ is represented by the identity map,
we get
$$
\pi _{0}(x)(e_{_{\Cal BP}})=\Id^{*}(x)=x,
$$
proving the lemma.

\enddemo

\sec
(5.2.9)
The product structure in absolute Hodge cohomology is given by a
morphism of complexes
$$
\Cal H^{*}(X,n,p)\otimes \Cal H^{*}(X,m,q)@>\cup >>
\Cal H^{*}(X,n+m,p+q),
$$
which induces a map of spaces
$$
\Cal K^{*}(X,n,p)\times \Cal K^{*}(X,m,q)@>\cup >>
\Cal K^{*}(X,n+m,p+q).
$$
The spaces $\Cal K(X,n,p)$ are naturally pointed by the element $0$.
Moreover $0\cup x=x\cup 0=0$. Therefore the
above map of spaces factors through:
$$
\Cal K^{*}(X,n,p)\times \Cal K^{*}(X,m,q)\longrightarrow
\Cal K^{*}(X,n,p)\land \Cal K^{*}(X,m,q)\longrightarrow
\Cal K^{*}(X,n+m,p+q).
$$

\proclaim{Lemma 5.12} Let $x\in H_{\Cal H}^{n}(\Cal BP,\Bbb R(p))$
and $y\in H_{\Cal H}^{m}(\Cal BP,\Bbb R(q))$. Then for any $i>0$ the map
$\pi _{i}(x\cup y)=0$.
\endproclaim
\demo{Proof} Let us write $E=\underline {\Hom} (N\frak U,\Cal BP)$.
Then the map $\pi (x\cup y)$ can be factored as
$$
\pi _{i}(E)@>\pi _{i}(\text{diag}) >>\pi _{i}(E\land E) @>>>
\pi _{i}(\Cal K(N\frak U,n,p)\land \Cal K(N\frak U,m,q))
@>>> \pi _{i}(\Cal K(N\frak U,u+m,p+q)).
$$
But since $S^{i}\land S^{i}=S^{2i}$ and for $i>0$, $\pi _{i}S^{2i}=0$,
the map $\pi _{i}(\text{diag})=0$.
\enddemo

\sec
(5.2.9) We are ready to prove theorem 5.2.
Let $i>0$ and $\sigma \in K_{i-1}(X,\frak U)$. By lemma 5.11, we have
that
$$
\align
\ch(\sigma )&=\pi_{i} (\ch (e_{_{\Cal BP}}))(\sigma ),\\
\rho (\sigma )&=\pi_{i} (\rho (e_{_{\Cal BP}}))(\sigma ).
\endalign
$$
Therefore
$$
\ch(\sigma )-\rho (\sigma )=
\pi_{i} (\ch (e_{_{\Cal BP}})-\rho (e_{_{\Cal BP}}))(\sigma ).
$$
By corollary 5.10, $\ch (e_{_{\Cal BP}})-\rho (e_{_{\Cal BP}})$ is
a sum of decomposable elements. Therefore by lemma 5.12.
$$
\ch(\sigma )=\rho (\sigma )
$$
concluding the proof of the theorem.
\enddemo

\heading {\S}6 Higher Bott-Chern forms.
\endheading

The higher Bott-Chern forms introduced in {\S}3 are differential forms
defined on $X\times (\Bbb P^{1})^{*}$. Nevertheless, the original
Bott-Chern forms (\cite {B-C}) and the higher Bott-Chern forms
introduced by Wang in \cite {Wan} are differential forms defined on $X$.
The aim of this section is to relate both notions of higher Bott-Chern
forms, in the case when $X$ is a proper smooth complex variety.
The main tool for this comparison will be an explicit quasi-isomorphism
$$
\widetilde {\frak H}^{*}(X,p)\longrightarrow \frak H^{*}(X,p).
$$
To this end we shall first introduce some complexes which compute
absolute Hodge homology and cohomology.

\sec
(6.1) Let us begin by introducing the complex where the simplest
Bott-Chern forms
are defined. This complex is a minor modification of the complex
used by Wang in \cite {Wan} (see also \cite {Bu 2}). The use of this
complex has been suggested by Deligne in \cite {De}. Let $X$ be a proper
smooth complex variety. We shall write
$$
E^{*}_{\Bbb R}(X)(p)=(2\pi i)^{p}E^{*}_{\Bbb R}(X).
$$

\sec
{\bf Definition 6.1.}
The complex $\frak W^{*}(X,p)$ is defined by
$$
\frak W^{n}(X,p)=
\cases
E^{n-1}_{\Bbb R}(X)(p-1)\cap
\dsize { \bigoplus \Sb p'+q'=n-1\\ p'<p,\ q'<p \endSb} E^{p',q'}(X),
&\qquad\text{ for } n\le 2p-1,\\
E^{n}_{\Bbb R}(X)(p)\cap
\dsize {\bigoplus \Sb p'+q'=n\\ p'\ge p,\ q'\ge p \endSb } E^{p',q'}(X)
\cap \Ker d,
&\qquad\text{ for } n= 2p,\\
0,
&\qquad\text{ for } n> 2p.\\
\endcases
$$
If $x\in \frak W^{n}(X,p)$ the differential $d_{\frak W}$ is given by
$$
d_{\frak W}x=\cases
-\pi (dx), & \text{for } n< 2p-1,\\
-2\partial \overline \partial x,\ &
\text{for } n=2p-1,\\
0, & \text{for } n= 2p,
\endcases
$$
where
$$\pi:E^{*}(X)\longrightarrow
E^{*}_{\Bbb R}(X)(p-1)\cap
\dsize { \bigoplus \Sb p'+q'=n-1\\ p'<p,\ q'<p \endSb} E^{p',q'}(X),
$$
is the projection.

\proclaim{Proposition 6.2} If $X$ is a proper smooth complex variety,
then
$$
H^{*}(\frak W^{*}(X,p))=H_{\Cal H}^{*}(X,\Bbb R(p)).
$$
\endproclaim
\demo{Proof}
Since $X$ is proper,
$$
H^{n}_{\Cal H}(X,\Bbb R(p))=
\cases
H^{n}_{\Cal D}(X,\Bbb R(p)),&\text{ for }n\le 2p,\\
0,&\text{ for }n>2p,
\endcases
$$
where $H^{n}_{\Cal D}(X,\Bbb R(p))$, denotes real Deligne cohomology of
$X$. Therefore the result follows from \cite {Bu 2 {\S}2}.
\enddemo

As in \cite {Bu 2}, we have morphisms of complexes
$$
\psi : \frak H^{*}(X,p)
\longrightarrow
\frak W^{*}(X,p)
$$
and
$$
\varphi :\frak W^{*}(X,p)
\longrightarrow
\frak H^{*}(X,p)
$$
given by
$$
\psi (a,f,\omega )=
\cases
\pi (\omega ),\ &\text{for }n\le 2p-1 \text{ and} \\
\sum_{i=p}^{n-p} a^{i,n-i}+\partial \omega
^{p-1,n-p+1}+(-1)^{p}\opartial \overline \omega ^{p-1,n-p+1},
\ &\text{for }n\ge 2p,
\endcases
$$
and
$$
\varphi (x)=
\cases
(\partial x^{p-1,n-p}-\overline \partial x^{n-p,p-1},2\partial
x^{p-1,n-p},x),\ &\text{for }n\le 2p-1 \text{ and} \\
(x,x,0),\ &\text{for } n\ge 2p,
\endcases
$$
where, if $x\in E_{X}^{*}$, then $x=\sum x^{p,q} $ is
the decomposition of $x$ in terms of pure type. The
morphisms
$\varphi $ and $\psi $ are homotopy equivalences inverse to each
other.

\sec
(6.2) In order to make the process of comparison clearer, we need an
auxiliary complex to compute
absolute Hodge cohomology, which is provided with a graded commutative
and associative product. It can be obtained by means of the Thom-Whitney
simple introduced by Navarro Aznar (see \cite {N 2} for the general
definition and properties of the Thom-Whitney simple).

Let $L^{*}_{1}$ be the differential graded commutative $\Bbb R$-algebra
of algebraic forms over $\Bbb A^{1}_{\Bbb R}$. Explicitly
$L^{0}_{1}=\Bbb
R[\epsilon ]$ and $L^{1}_{1}=\Bbb R[\epsilon ]d\epsilon $. Let $\delta
_{0}:L^{*}_{1}\longrightarrow \Bbb R$ (resp. $\delta _{1}$) be the
evaluation at $0$ morphism (resp. evaluation at $1$).

\sec
{\bf Definition 6.3.} Let $X$ be a smooth complex variety. The
Thom-Whitney simple of the absolute Hodge complex, denoted by
$\frak H^{*}_{TW}(X,p)$, is the subcomplex of
$$
\left((2\pi i)^{p}\widehat W_{2p}E_{\log}^{*}(X)_{\Bbb R}\oplus
\widehat W_{2p} \cap F^{p}\, E^{*}_{\log}(X)\oplus
(L^{*}_{1}\underset \Bbb R \to \otimes \widehat W_{2p} E_{\log}^{*}(X))
\right)
$$
formed by the elements $(r,f,\omega)$ such that
$$
\align
\omega (0)&=(\delta _{0}\otimes \Id)\, (\omega )=r,\\
\omega (1)&=(\delta _{1}\otimes \Id)\, (\omega )=f.
\endalign
$$

Let $E$ and $I$ be the morphisms of complexes
$$
\frak H^{*}_{TW}(X,p)
\matrix  I\\ \longrightarrow \\ \longleftarrow \\E \endmatrix
\frak H^{*}(X,p)
$$
given by
$$
\align
E(r,f,\omega )&=(r,f,\epsilon \otimes f+(1-\epsilon )
\otimes r+d\epsilon \otimes \omega ),
\\
I(r,f,\omega )&=(r,f,\int_{0}^{1}\omega  ),
\endalign
$$
where the integration symbol means formal integration with
respect to the variable $\epsilon $.
These morphisms are homotopy equivalences (see \cite {N 2}).

We shall denote by $I'$ the composition
$$
\frak H_{TW}^{*}(X,*)@>I>> \frak H^{*}(X,*)
@>\psi >> \frak W^{*}(X,*),
$$
and by $E'$ the composition
$$
\frak W^{*}(X,*)@>\varphi >> \frak H^{*}(X,*)
@>E >> \frak H_{TW}^{*}(X,*).
$$
The morphisms $I'$ and $E'$ are also homotopy equivalences inverse to
each other.

We can define a product
$$
\frak H^{n}_{TW}(X,p)
\otimes
\frak H^{m}_{TW}(X,q)
@>\cup>>
\frak H^{n+m}_{TW}(X,p+q),
$$
by
$$
(r,f,\omega)
\cup
(r',f',\omega')=
(r\land r',f\land f',\omega \land \omega ').
$$
This product is associative, graded commutative and satisfies the
Leibnitz rule. Therefore
$$
\frak H^{*}_{TW}(X,*)=\bigoplus_{p}
\frak H^{*}_{TW}(X,p)
$$
is a differential associative graded commutative algebra. Moreover, the
$\Bbb R$-algebra structure induced in $H^{*}_{\Cal H}(X,\Bbb R(p))$ by
this product coincides with the $\Bbb R$-algebra structure introduced by
Beilinson (\cite {Be}).

\sec
(6.3) Let us give the homology analogue of the last complex. This is
done by means of currents. For a proper smooth complex variety $X$, let
$D_{*.*}(X)$ be the double chain complex of complex currents over $X$,
let $D_{*}(X)$ be the associated single complex, and let $D^{\Bbb
R}_{*}(X)$ be the real subcomplex. We shall write
$$
F_{p}D_{*}(X)=\bigoplus _{p'\le p}D_{p',*}.
$$
Let $\tau ^{\ge 2p}D_{*}(X)$ be the subcomplex
$$
\tau ^{\ge 2p}D_{n}(X)=\cases
D_{n}(X), &\text{ if } n>2p,\\
\Ker (d), &\text{ if } n=2p,\\
0, &\text{ if } n<2p.
\endcases
$$
Since $X$ is proper, the filtration $\tau $ plays the role of the
d\'ecal\'ee weight filtration.

Let $L_{*}^{1}$ be the chain complex defined by $L_{k}^{1}=L^{-k}_{1}$
(see 6.2). We shall denote by $\delta _{0}$ and $\delta _{1}$ the
evaluation at $0$ and $1$ as in (6.2).

\sec
{\bf Definition 6.4.} Let
$\frak H_{*}^{TW}(X,p)$ be the subcomplex of
$$
\left((2\pi i)^{-p}\tau ^{\ge 2p}D_{*}^{\Bbb R}(X)\oplus
\tau ^{\ge 2p} \cap F_{p}\, D_{*}(X)\oplus
(L_{*}^{1}\underset \Bbb R \to \otimes \tau ^{\ge 2p} D_{*}(X))
\right)
$$
formed by the elements $(r,f,\omega)$ such that
$$
\align
\omega (0)&=(\delta _{0}\otimes \Id)\, (\omega )=r,\\
\omega (1)&=(\delta _{1}\otimes \Id)\, (\omega )=f.
\endalign
$$

The homology of the complex $\frak H_{*}^{TW}(X,p)$ is the absolute
Hodge homology of $X$.

\sec
(6.4) The last complex we introduce is an analogue of
$\widetilde {\frak H}^{*}(X,p)$, replacing $\frak H^{*}(X,p)$ by
$\frak H_{TW}^{*}(X,p)$.
We shall denote by
$
\frak H_{\Bbb P,TW}^{*,*}(X,p)
$
the double complex given by
$$
\frak H_{\Bbb P,TW}^{r,n}(X,p)=
\frak H^{r}_{TW}(X\times (\Bbb P^{1})^{-n},p),
$$
with differentials $$
\align
d'&=d_{\frak H},\\
d''&=\sum (-1)^{i+j} d^{j}_{i}.
\endalign
$$
Then the double complex
$\widetilde {\frak H}_{TW}^{*,*}(X,p)$ is given by
$$
\widetilde {\frak H}_{TW}^{r,n}(X,p)=
\frak H^{r,n}_{\Bbb P,TW}(X,p)\left /
\sum_{i=1}^{n}
s_{i}\left(\frak H^{r,n+1}_{\Bbb
P,TW}(X,p)\right)\,\oplus\,
\omega _{i}\land s_{i}\left(\frak H^{r-2,n+1}_{\Bbb P,TW}(X,p-1)\right).
\right.
$$
Finaly let
$\widetilde {\frak H}_{TW}^{*}(X,p)$ be the associated simple complex.
The differential of this complex will be denoted by $d$.

Observe that the homotopy equivalences $E$ and $I$ induce homotopy
equivalences $$
\widetilde {\frak H}^{*}_{TW}(X,p)
\matrix  I\\ \longrightarrow \\ \longleftarrow \\E \endmatrix
\widetilde {\frak H}^{*}(X,p).
$$

\sec
(6.5) In order to pull forms in $X\times (\Bbb P^{1})^{n}$ down to $X$,
we need some differential forms on $X\times (\Bbb
P^{1})^{n}$ which will play a role similar to the currents ''integration
along the standard simplex''.

Let $(x:y)$ be homogeneous coordinates of $\Bbb
P^{1}$, and let $t=x/y$ be the absolute coordinate of $\Bbb
P^{1}$. Let us write $\Bbb C^{*}=\Bbb P^{1}_{\Bbb
C}-\{0,\infty \}$. Let
$$
\align
\lambda &=\frac{1}{2}E'(\log t\overline t)\\
&=\frac{1}{2} \left(\frac{dt}{t}-\frac{d\overline t}{\overline t},
2\frac{dt}{t},
(\varepsilon  +1)\otimes \frac{dt}{t}+
(\varepsilon -1) \otimes \frac{d\overline t}{\overline t}
+d\varepsilon \otimes \log t\overline t
\right)\\
&\in \frak H^{1}_{TW}(\Bbb C^{*},1).
\endalign
$$

Let us consider the open subset
$(\Bbb C^{*})^{n}\subset X\times (\Bbb P^{1})^{n}$. Let us denote by
by $p_{i}:(\Bbb C^{*})^{n}\longrightarrow \Bbb C^{*}$, $i=1,\dots ,n$
the projections over the $i$-th factor.
Let us write $\lambda _{i}=p_{i}^{*}\lambda $.

\sec
{\bf Definition 6.5.} Let $W_{n}\in \frak H^{n}_{TW}((\Bbb
C^{*})^{n},n)$ be the form defined by
$$
W_{n}=\lambda _{1}\cup \dots \cup \lambda _{n}.
$$

\sec
(6.6) Since the forms $W_{n}$ will play a central role, let us present
a more explicit description. Let us write
$W_{n}=(W_{n}^{1},W_{n}^{2},W_{n}^{3})$.
Then
$$
\aligned
W_{n}^{1}&=\frac{1}{2^{n}} \bigwedge_{i=1}^{n}
\left(
\frac{dt_{i}}{t_{i}}-\frac{d\overline t_{i}}{\overline t_{i}}\right)\\
W_{n}^{2}&=\bigwedge_{i=1}^{n}\frac{dt_{i}}{t_{i}}\\
W_{n}^{3}&=\frac{1}{2^{n}}\bigwedge_{1=1}^{n} \left(
(\epsilon +1)\otimes \frac{dt_{i}}{t_{i}}+
(\epsilon-1) \otimes \frac{d\overline t_{i}}{\overline t_{i}}
+d\epsilon \otimes \log t_{i}\overline t_{i}\right).
\endaligned
$$
Let $\frak S_{n}$ denote the symmetric group. Let us write, for
$i=0,\dots ,n$,
$$
P^{i}_{n}=\sum_{\sigma \in \frak S_{n}}(-1)^{\sigma }
\frac{dt_{\sigma (1)}}{t_{\sigma (1)}}
\land \dots \land\frac{dt_{\sigma (i)}}{t_{\sigma (i)}}\land
\frac{d\overline t_{\sigma (i+1)}}{\overline t_{\sigma (i+1)}}
\land \dots \land
\frac{d\overline t_{\sigma (n)}}{\overline t_{\sigma (n)}},
$$
and, for $i=1,\dots ,n$,
$$
S^{i}_{n}=\sum_{\sigma \in \frak S_{n}}(-1)^{\sigma }
\log (t_{\sigma (1)} \overline t_{\sigma (1)})
\frac{dt_{\sigma (2)}}{t_{\sigma (2)}}\land \dots
\land\frac{dt_{\sigma (i)}}{t_{\sigma (i)}}\land
\frac{d\overline t_{\sigma (i+1)}}{\overline t_{\sigma (i+1)}}
\land \dots \land
\frac{d\overline t_{\sigma (n)}}{\overline t_{\sigma (n)}}.
$$
Then we have
$$
\aligned
W_{n}^{1}&=\frac{1}{2^{n}} \sum_{i=0}^{n}(-1)^{n-i}\frac{1}{i!(n-i)!}
P^{i}_{n},\\
W_{n}^{3}&=\frac{1}{2^{n}}
\sum_{i=0}^{n}
\frac{(\epsilon+1)^{i}(\epsilon-1)^{n-i}}{i!(n-i)!}\otimes P^{i}_{n}+
\frac{1}{2^{n}} \sum_{i=1}^{n}
\frac{(\epsilon+1)^{i-1}(\epsilon-1)^{n-i}}{(i-1)!(n-i)!}d\epsilon
\otimes S^{i}_{n}.
\\
\endaligned
$$

\sec
(6.7) We are not as interested in the forms $W_{n}$, as in their
associated currents. Let
$\omega \in E^{r}_{(\Bbb P^{1})^{n}}$. Let us denote by $[\omega ]\in
D_{2n-r}((\Bbb P^{1})^{n})$ the current defined by
$$
[\omega ](\varphi )=\frac{1}{(2\pi i)^{n}}
\int_{(\Bbb P^{1})^{n}} \varphi \land \omega .
$$
If $a\otimes \omega \in L^{*}_{1}\otimes E^{r}_{(\Bbb P^{1})^{n}}$
we write
$$
\align
[a\otimes \omega ]&=a\otimes [\omega ]\\
&\in L^{1}_{*}\otimes D_{2n-r}((\Bbb P^{1})^{n}).
\endalign
$$
In this way we obtain a map
$$
\matrix
\frak H_{TW}^{r}((\Bbb P^{1})^{n},p)&
\longrightarrow&
\frak H^{TW}_{2n-r}((\Bbb P^{1})^{n},n-p)\\
(r,f,\omega )&\longmapsto&([r],[f],[\omega ])
\endmatrix
$$

This definition can be extended to any locally integrable differential
form.

\sec
{\bf Definition 6.6.}
We shall denote by $[W_{n}]$ the element of
$ \frak H_{n}^{TW}((\Bbb P^{1})^{n},0) $ given by
$$
[W_{n}]=([W_{n}^{1}],[W_{n}^{2}],[W_{n}^{3}]).
$$

The following result exhibits the analogy between the currents
``integration along the
standard simplex'' and the currents $[W_{n}]$.

\proclaim{Proposition 6.7} The currents $[W_{n}]$ satisfy the relation
$$
d[W_{n}]=\sum_{i=1}^{n}\sum_{j=0,1}(-1)^{i+j} (d^{i}_{j})_{*}[W_{n-1}].
$$
\endproclaim
\demo{Proof} Formally this proposition is the Leibnitz rule.
To prove it we can work component by component.
By a standard residue argument:
$$
\align
d[W_{n}^{2}]&=d\left [\bigwedge_{i=1}^{n}\frac{dt_{i}}{t_{i}}\right]\\
&=\sum_{i=1}^{n}\sum_{j=0}^{1}(-1)^{i+j}(d^{i}_{j})_{*}[W_{n-1}^{2}].
\endalign
$$
By the same argument and taking some care with permutations one
sees
$$
\align
d[P^{i}_{n}]&=
\sum_{k=1}^{n}\sum_{j=0}^{1}(-1)^{k+j}
(d^{k}_{j})_{*}(i[P^{i-1}_{n-1}]-(n-i)[P^{i}_{n-1}]),\\
d[S^{i}_{n}]&=[P^{i}_{n}]+[P^{i-1}_{n}]+
\sum_{k=1}^{n}\sum_{j=0}^{1}(-1)^{k+j}
(d^{k}_{j})_{*}((i-1)[S^{i-1}_{n-1}]-(n-i)[S^{i}_{n-1}]).\\
\endalign
$$
The proposition follows from the above formulas and the explicit
description of $W_{n}^{1}$ and $W_{n}^{3}$ given in (6.6).
\enddemo

\sec
(6.8) Acting component by
component, the currents $[W_{n}]$ induce morphisms
$$
[W_{n}]: \frak H_{\Bbb P,TW}^{r,-n}(X,p)=
\frak H_{TW}^{r}(X\times (\Bbb P^{1})^{n},p)
\longrightarrow
\frak H_{TW}^{r-n}(X,p).
$$

\proclaim{Lemma 6.8} The morphisms $[W_{n}]$ factorize through morphisms
$$
[W_{n}]:\widetilde {\frak H}^{r,-n}_{TW}(X,p)
\longrightarrow
\frak H^{r-n}_{TW}(X,p).
$$
\endproclaim
\demo{Proof} Let us denote by $\sigma _{i}$ the automorphism of
$(\Bbb P^{1})^{n}$ given by
$$
\sigma _{i}((x_{1}:y_{1}),\dots ,(x_{i}:y_{i}),\dots (x_{n}:y_{n}))
=
((x_{1}:y_{1}),\dots ,(y_{i}:x_{i}),\dots (x_{n}:y_{n})).
$$
Then  $(\sigma _{i})_{*}[W_{n}]=-[W_{n}]$, for $i=1,\dots ,n$. On the
other hand, if
$$
\eta \in
s_{i}\left(\frak H^{r,n+1}_{\Bbb
P,TW}(X,p)\right)\,\oplus\,
\omega _{i}\land s_{i}\left(\frak H^{r-2,n+1}_{\Bbb P,TW}(X,p-1)\right)
$$
then $(\sigma _{i})^{*}\eta=\eta$. Therefore
$$
[W_{n}]\eta=-(\sigma _{i})_{*}[W_{n}]\eta=
-[W_{n}](\sigma _{i})^{*}\eta=-[W_{n}]\eta.
$$
Hence $[W_{n}]\eta=0$ proving the result.
\enddemo

\sec
{\bf Definition 6.9.} Let $W_{TW}$ be the morphism
$$
W_{TW} :\widetilde {\frak H}_{TW}^{*}(X,p)\longrightarrow
\frak H^{*}_{TW}(X,p)
$$
given, for $\eta \in \widetilde {\frak H}_{TW}^{r,-n}(X,p)$, by
$$
W_{TW} (\eta)=\cases [W_{n}]\eta, & \text{ if }n>0,\\
\eta, & \text{ if }n=0.\endcases
$$

\proclaim{Proposition 6.10} The morphism $W_{TW} $ is a morphism of
complexes. Moreover it is a quasi-isomorphism.
\endproclaim
\demo{Proof} The fact that is a morphism of complexes is a consequence
of proposition 6.7.
Let
$$\iota':
\frak H^{*}_{TW}(X,p)\longrightarrow
\widetilde {\frak H}_{TW}^{*}(X,p)$$
be the morphism induced by the
equality $\frak H^{*}_{TW}(X,p) =\widetilde {\frak H}_{TW}^{*,0}(X,p)$.
We have that $\iota= I\circ \iota' \circ E$, where $\iota$ is the
quasi-isomorphism defined in proposition 1.2. Therefore
$\iota '$ is a
quasi-isomorphism. Since $W_{TW} \circ \iota '=\Id$ we have that $W_{TW}
$ is also a quasi-isomorphism.
\enddemo

\sec
{\bf Definition. 6.11.} Let us denote by $W$ the morphism
$$
        W= I  \circ    W_{TW}  \circ  E:
\widetilde {\frak H}^{*}(X,p)\longrightarrow \frak H^{*}(X,p).
$$

Observe that $W$ is also a quasi-isomorphism. Summarizing, we have the
following diagram of complexes and quasi-isomorphisms.
$$
\CD
\widetilde {\frak H}^{*}(X,p)@>W>> \frak H^{*}(X,p)@>\psi >> \frak
W^{*}(X,p) \\
@V E VV  @A I AA @.\\
\widetilde {\frak H}_{TW}^{*}(X,p)@>W_{TW}>> \frak H^{*}_{TW}(X,p)@.
\endCD
$$

\sec
(6.9) The diagram above allows us to define different kinds of higher
Bott-Chern forms. For instance let us recover
the original
definition of higher Bott-Chern forms due to Wang (\cite {Wan})
and
the classical Bott-Chern
forms.

\sec
{\bf Definition 6.12.} Let $\overline {\Cal F}$ be an exact
metrized $n$-cube. We shall also call the Bott-Chern form of $\overline
{\Cal F}$ the form
$$
\chw{n}(\overline {\Cal F})_{W}=\psi \circ W(\chw{n}(\overline {\Cal
F})_{\Cal H}). $$

One may compute these forms directly using the following
result.

\proclaim{Proposition 6.13} Let $\overline {\Cal F}$ be an
emi-$n$-cube. Then
$$
\chw{n}(\overline {\Cal F})_{W}=\frac{1}{(2\pi i)^{n}}\int_{(\Bbb
P^{1})^{n}} \chw{0}(\tr_{n}(\overline {\Cal F}))\land I'(W_{n}).
$$
\endproclaim
\demo{Proof} This result is consequence of the following facts
\roster
\item The morphism $I'$ is functorial.
\item For any smooth complex variety $Z$, if $\omega \in \frak
H_{TW}^{2p}(Z,p)$ and $\eta\in \frak H_{TW}^{*}(Z^{*},*)$, then
$I'(\omega \cup \eta)=I'(\omega )\land I'(\eta)$ (see \cite {Bu 2}).
\item $I'\circ E'=\Id$. Therefore
$$
I'(E(\chw{0}(\tr_{n}(\overline {\Cal F}))_{\Cal H}))=
I'(E'(\chw{0}(\tr_{n}(\overline {\Cal F}))))=
\chw{0}(\tr_{n}(\overline {\Cal F})).
$$
\endroster
\enddemo

Up to a normalization factor, the formula given in proposition 6.13
is the original definition due to Wang (\cite {Wan}).
To see this, let us compute explicitly
$I'(W_{n})\in \frak W^{n}((\Bbb C^{*})^{n},n)$.

\proclaim{Proposition 6.14}
$$
I'(W_{n})=\frac{(-1)^{n}}{2\,n!}\sum_{i=1}^{n}(-1)^{i-1}S^{i}_{n}.
$$
\endproclaim
\demo{Proof}
Since $W_{n}\in
\frak H^{n}_{TW}((\Bbb C^{*})^{n},n)$, by (6.1) and (6.2), we have
$$
I'(W_{n})=\pi \left(\int_{0}^{1}W_{n}^{3}\right),
$$
where the integral symbol means integration with respect to the variable
$\epsilon$, and $\pi $ is the projection
$$
\pi :E^{n-1}_{(\Bbb C^{*})^{n}}\longrightarrow
(2\pi i)^{n-1}E^{n-1}_{(\Bbb C^{*})^{n},\Bbb R}.
$$
This projection is given by $\pi (z)=(z+(-1)^{n-1}\overline z)/2$.
Therefore
$$
I'(W_{n})=\frac{1}{2^{n+1}}\sum_{i=1}^{n}
\int_{0}^{1}
\frac{(\epsilon+1)^{i-1}(\epsilon-1)^{n-i}}{(i-1)!(n-i)!}d\epsilon\,
\left(S^{i}_{n}+(-1)^{n-1} \overline S^{i}_{n}
\right).
$$
But $\overline S^{i}_{n}=S^{n-i+1}_{n}$. Then, joining the terms with
$S^{i}_{n}$, and taking into account that
$$
(-1)^{n-1}\int_{0}^{1}
\frac{(\epsilon+1)^{n-i}(\epsilon-1)^{i-1}}{(i-1)!(n-i)!}d\epsilon=
\int_{-1}^{0}
\frac{(\epsilon+1)^{i-1}(\epsilon-1)^{n-i}}{(i-1)!(n-i)!}d\epsilon,
$$
we have that
$$
I'(W_{n})=\frac{1}{2^{n+1}}\sum_{i=1}^{n}\left(
\int_{-1}^{1}
\frac{(\epsilon+1)^{i-1}(\epsilon-1)^{n-i}}{(i-1)!(n-i)!}d\epsilon\,
S^{i}_{n}
\right).
$$
But
$$
\int_{-1}^{1}
\frac{(\epsilon+1)^{i-1}(\epsilon-1)^{n-i}}{(i-1)!(n-i)!}d\epsilon=
\frac{(-1)^{n+i-1}2^{n}}{n!},
$$
proving the result.
\enddemo

The following result is a direct
consequence of the definitions.
\proclaim{Proposition 6.15} Let $X$ be a proper smooth complex variety.
Let
$$
\overline {\xi }:
0@>>>\overline {\Cal F}
@>>>\overline {\Cal G}
@>>>\overline {\Cal H}
@>>>0
$$
be an exact sequence of locally free sheaves over $X$. Let us denote by
$\bc(\overline {\xi
})$ the Bott-Chern form of $\overline {\xi }$ as defined by Bismut,
Gillet and Soul\' e (\cite {B-G-S}, \cite {G-S 1}). Then
$$
\chw{1}(\overline {\xi })_{W}=-\frac{1}{2} \bc \overline {\xi }\qquad
\mod(\Im \partial +\Im \overline \partial ).
$$

\endproclaim

\sec (6.10) The use of the Thom-Whitney simple for absolute Hodge
cohomology, besides giving a way to construct the currents $W_{n}$,
allows us to define a multiplicative theory of Bott-Chern forms.

\sec
{\bf Definition 6.16.} Let $\overline {\Cal F}$ be an exact metrized
$n$-cube. We shall call the multiplicative Bott-Chern form of $\overline
{\Cal F}$ the form
$$
\chw{n}(\overline {\Cal F})_{TW}=W_{TW}(E(\chw{n}(\overline
{\Cal F})_{\Cal H})).
$$
In particular, if $\overline {\Cal F}$ is a hermitian locally free
sheaf, then
$$
\chw{0}(\overline {\Cal F})_{TW}=
E(\chw{0}(\overline {\Cal F})_{\Cal H}).
$$
On the other hand, if $\overline {\Cal F}$ is an emi-$n$-cube, then
$$
\chw{n}(\overline {\Cal F})_{TW}
=\frac{1}{(2\pi i)^{n}}\int_{(\Bbb P^{1})^{n}}
\chw{0}(\tr_{n}(\overline {\Cal F}))_{TW}\cup W_{n}.
$$

\sec
{\bf Definition 6.17.} Let $\overline {\Cal F}$ be a metrized exact
$n$-cube and let $\overline {\Cal G}$ be a metrized exact
$m$-cube. Then $\overline {\Cal F}\otimes \overline {\Cal G}$ is the
metrized exact $n+m$-cube given by
$$
(\overline {\Cal F}\otimes \overline {\Cal G})_{i_{1},\dots ,i_{n+m}}=
(\overline {\Cal F})_{i_{1},\dots ,i_{n}}\otimes
(\overline {\Cal G})_{i_{n+1},\dots ,i_{n+m}}
$$
with the obvious morphisms and metrics.

\proclaim{Proposition 6.18}
Let
$\overline {\Cal F}$ (resp.  $\overline {\Cal G}$) be a metrized exact
$n$-cube (resp. $m$-cube). Then
$$
\chw{n+m}(\overline {\Cal F}\otimes \overline {\Cal G})_{TW}
=
\chw{n}(\overline {\Cal F})_{TW}\cup\chw{m}(\overline {\Cal G})_{TW}.
$$
\endproclaim
\demo{Proof} We may assume that
$\overline {\Cal F}$ and $\overline {\Cal G}$ are emi-cubes.

Let $\pi _{1}:X\times (\Bbb P^{1})^{n+m}\longrightarrow
X\times (\Bbb P^{1})^{n}$ be the projection over the first
$n$-projective lines and let
$\pi _{2}:X\times (\Bbb P^{1})^{n+m}\longrightarrow
X\times (\Bbb P^{1})^{m}$ be the projection over the last
$m$-projective lines.

\proclaim{Lemma 6.19}
Let
$\overline {\Cal F}$ (resp.  $\overline {\Cal G}$) be an
emi-$n$-cube (resp. emi-$m$-cube). Then
$$
\tr_{n+m}(\overline {\Cal F}\otimes \overline {\Cal G})
=
\pi _{1}^{*}\tr_{n}(\overline {\Cal F})\otimes \pi
_{2}^{*}\tr_{m}(\overline {\Cal G}).
$$
\endproclaim
\demo{Proof} By {\S}3, (3.7), it is enough to show that, if $m\ge 1$,
then $$
\tr_{1}(\overline {\Cal F}\otimes \overline {\Cal G})
=
\overline {\Cal F}\otimes \tr_{1}(\overline {\Cal G}),
$$
and if $m=0$, then
$$
\tr_{1}(\overline {\Cal F}\otimes \overline {\Cal G})
=
\tr_{1}(\overline {\Cal F})\otimes \overline {\Cal G}.
$$
Since $\tr_{1}$ is computed in each edge separately, it is enough to
prove the case $n=1$, $m=0$, but this case follows directly from the
definition.
\enddemo

Using lemma 6.19, the multiplicativity and functoriality of the Chern
character forms and the definition of the forms $W_{n}$, we have:
$$
\align
\chw{n+m}(&\overline {\Cal F}\otimes \overline {\Cal G})_{TW}=\\
&=\frac{1}{(2\pi i)^{n+m}}\int_{(\Bbb P^{1})^{n+m}}
\chw{0}(\pi _{1}^{*}\tr_{n}(\overline {\Cal F})\otimes
\pi _{2}^{*}\tr_{m}(\overline {\Cal G}))_{TW}\cup W_{n+m}\\
&=\frac{1}{(2\pi i)^{n+m}}\int_{(\Bbb P^{1})^{n+m}}
\pi _{1}^{*}\chw{0}(\tr_{n}(\overline {\Cal F}))_{TW}\cup
\pi _{2}^{*}\chw{0}(\tr_{m}(\overline {\Cal G}))_{TW}\cup
\pi _{1}^{*}W_{n}\cup \pi _{2}^{*}W_{m}\\
&=
\chw{n}(\overline {\Cal F})_{TW}\cup
\chw{m}(\overline {\Cal G})_{TW}.
\endalign
$$
\enddemo

\heading References.
\endheading
\widestnumber\key {G-N-P-PP}

\ref\key Be \by Beilinson, A. A. \paper Notes on absolute Hodge
cohomology
\jour Contemporary Mathematics\vol 55  I\yr
1986 \pages 35--68\endref

\ref\key B-G-S \by Bismut, J.M., Gillet, H. and Soul\'e, C. \paper
Analytic torsion and holomorphic determinant bundles I: Bott-Chern forms
and analytic torsion,
\jour Comm. Math. Phys. \vol 115 \yr 1988 \pages 49--78\endref

\ref\key B-C \by Bott, R. and Chern, S. S. \paper Hermitian vector
bundles and the equidistribution of zeroes of their holomorphic sections
\jour Acta Math \vol 114 \yr 1968 \pages
71--112\endref

\ref\key B-K \by Bousfield A.K. and Kan D.M.
\book Homotopy limits, completions and localizations.
LNM 304
\publ Springer-Verlag\yr 1972\endref

\ref\key Bu 1 \by Burgos J. I.\paper
A C$^{\infty }$ logarithmic Dolbeault complex \jour Compos. Math.
\vol 92 \yr 1994 \pages 61--86 \endref

\ref\key Bu 2 \by Burgos J. I.\paper
Arithmetic Chow rings and Deligne Beilinson cohomology\jour J. of
Algebraic Geometry \toappear \endref

\ref\key Ch \by Chern S. S.\paper
Characteristic classes of hermitian manifolds \jour Ann. of Math.
\vol 47 \yr 1946 \pages 85--121 \endref

\ref\key De \by Deligne P.\paper
Le d\'eterminant de la cohomologie\jour Contemporary Mathematics
\vol 67 \yr 1987 \pages 93--177 \endref

\ref\key Gi \by Gillet, H. \paper
Riemann-Roch theorems for higher algebraic K-theory
\jour Advances in Math.\vol 40 \yr 1981 \pages 203--289\endref

\ref\key G-S 1 \by Gillet, H. and Soul\'e, C. \paper Characteristic
classes for algebraic vector bundles with hermitian metric I and II
\jour Annals of Math. \vol 131 \yr 1990 \pages 163--238\endref

\ref\key G-S 2\by Gillet, H. and Soul\'e, C. \paper Arithmetic
intersection
theory \jour Publ. Math. IHES \vol 72 \yr 1990 \pages 93--174\endref

\ref\key Gr \by Griffiths, P. A. \paper On a theorem of Chern
\jour Illinois J. Math. \vol 6 \yr 1962 \pages 468--479\endref

\ref\key G-D \by Grothendieck, A. and Dieudonn\' e, J.A.
\book El\'ements de g\'eom\'etrie alg\'e\-bri\-que
\publ Sprin\-ger-Verlag\yr 1971\endref

\ref\key GNPP \by Guillen, F., Navarro Aznar, V., Pascual-Gainza,
P. and Puerta, F. \book Hyperr\'esolu\-tions cubiques et descente
cohomologique. LNM 1335
\publ Sprin\-ger-Verlag\yr 1988\endref

\ref \key Ka 1 \by Karoubi, M.\paper Homologie cyclique et $K$-th\'eorie
\jour Ast\'erisque
\vol 149\yr 1987 \endref

\ref \key Ka 2 \by Karoubi, Th\'eorie g\'en\'erale des classes
caract\'eristiques secondaries.
\jour $K$-Theory
\vol 4\yr 1990 \pages 55--88\endref

\ref \key Le \by Lecomte, F. \paper Op\'erations d'Adams en
$K$-Th\'eorie alg\'ebrique\jour Strasbourg Preprint
\yr 1996\endref

\ref \key Lo 1 \by Loday, J.L. \paper Homotopie des spaces de
concordances (d'apr\`es F. Waldhausen)
\inbook Sem. Bourbaki, 1977/78 N. 516,  LNM
710\pages 187--205\publ Springer-Verlag \yr 1979\endref

\ref \key Lo 2 \by Loday, J.L. \paper Spaces with finitely many
non-trivial homotopy groups
\jour J. of Pure and Applied Algebra
\vol 24\yr 1982 \pages 179--202\endref

\ref\key Mas \by Massey, W. S.
\book Singular Homology Theory, GTM 70
\publ Springer Verlag \yr 1980\endref

\ref\key M-S \by Milnor, J. W., Stasheff, J. D.
\book Characteristic Classes, Annals of Mathematics Studies 76
\publ Princeton University Press\yr 1974\endref

\ref \key N 1\by Navarro Aznar, V. \paper Sur les multiplicit\'es de
Schubert locales des faisceaux alg\'ebri\-ques coh\'erents
\jour Comp. Math.
\vol 48\yr 1983 \pages 311--326\endref

\ref\key N 2\by Navarro Aznar, V. \paper Sur la th\'eorie de
Hodge-Deligne\jour Invent. Math.\vol 90 \yr 1987 \pages 11--76\endref

\ref \key Ri \by Riemenschneider, O. \paper Characterizing Moishezon
spaces by almost positive coherent analytic sheaves
\jour Math. Z.
\vol 123\yr 1971 \pages 263--284\endref

\ref \key Ro \by Rossi, H. \paper Picard variety of an isolated singular
point
\jour Rice Univ. Studies
\vol 54\yr 1968 \pages 63--73\endref

\ref \key Sch \by Schechtman \paper On the delooping of Chern character
and Adams operations \inbook K-theory, arithmetic and geometry. LNM
1289\pages 265--319\publ Springer-Verlag \yr 1987\endref

\ref\key So \by Soul\'e et als.
\book Lectures on Arakelov Geometry, Cambridge studies in advanced
mathematics 33\publ Cambridge University Press\yr 1992\endref

\ref \key Th \by Thomason, R.
\paper Algebraic K-theory and \'etal\'e cohomology
\jour Ann. Sci. Ec. Norm. Sup.
\vol 18\yr 1985 \pages 437--552\endref

\ref \key Wal \by Waldhausen, F. \paper Algebraic K-theory of
generalized free products
\jour Ann. Math.
\vol 108\yr 1978 \pages 135--256\endref

\ref\key Wan \by Wang X. \book  Higher-order characteristic classes in
arithmetic geometry \bookinfo Thesis Harvard \yr 1992\endref

\enddocument